\documentclass[11pt]{smfart}

\usepackage{amssymb}
\usepackage{amsmath, amsthm, amsopn, amsfonts}
\usepackage[dvips]{graphicx}
\usepackage{graphpap}
\usepackage[francais, english]{babel}

\def\11{{\rm 1~\hspace{-1.4ex}l} }
\def\R{\mathbb R}
\def\C{\mathbb C}
\def\Z{\mathbb Z}
\def\N{\mathbb N}

\def\T{\mathbb T}
\def\O{\mathbb O}

\def\cal{\mathcal}

\newtheorem{theoreme}{Theorem}

\newtheorem*{lemmecite}{Lemma}
\renewcommand{\theequation}{\arabic{equation}}

\newtheorem{proposition}{Proposition}
\newtheorem{lemme}[proposition]{Lemma}

\newtheorem{remarque}[proposition]{Remark}

\numberwithin{equation}{section}
\numberwithin{proposition}{section}
\newcommand{\rien}[1]{\relax}

\def\smfbyname{by}
\begin{document}
\selectlanguage{english}
\title[Nonlinear Schr\"odinger equations on 3-manifolds]
{Multilinear Eigenfunction Estimates And Global Existence For The Three Dimensional Nonlinear  Schr\"Odinger Equations }

\author{N. Burq}
\address{Universit{\'e} Paris Sud,
Math{\'e}matiques,
B{\^a}t 425, 91405 Orsay Cedex et Institut Universitaire de France}
\email{Nicolas.burq@math.u-psud.fr}
\urladdr{http://www.math.u-psud.fr/~burq}
\author{P. G{\'e}rard}
\address{Universit{\'e} Paris Sud,
Math{\'e}matiques,
B{\^a}t 425, 91405 Orsay Cedex}
\email{Patrick.gerard@math.u-psud.fr}
\author{N. Tzvetkov}
\address{Universit{\'e} Paris Sud,
Math{\'e}matiques,
B{\^a}t 425, 91405 Orsay Cedex}
\email{Nikolay.tzvetkov@math.u-psud.fr}
\urladdr{http://www.math.u-psud.fr/~tzvetkov}
\selectlanguage{english}
\begin{abstract}
We study nonlinear  Schr\"odinger equations, posed on a three dimensional Riemannian manifold $M$.
We prove global existence of strong $H^1$ solutions on $M=S^3$ and $M=S^2\times S^1$ as far as the nonlinearity
is defocusing and sub-quintic and thus we extend results of Ginibre-Velo  and Bourgain  who treated
the cases of the Euclidean space $\R^3$ and the torus $\T^3= \mathbb{R}^3/\mathbb{Z}^3$ respectively. The main ingredient in our argument is
a new set of multilinear estimates for spherical harmonics. 
\end{abstract} 
\begin{altabstract}
On \'etudie l'\'equation de Schr\"odinger sur une vari\'et\'e de dimension trois $M$. On d\'emontre l'existence globale en temps de solutions fortes $H^1$ si $M= S^3$ ou $S^2 \times S^1$, pour les non lin\'earit\'es sous quintiques et d\'efocalisantes. On \'etend ainsi les r\'esultats de Ginibre et Velo et Bourgain qui ont trait\'e les cas de l'espace euclidien $R^3$ et du tore $\mathbb{T}^3= \mathbb{R}^3/\mathbb{Z}^3$ respectivement. L'ingr\'edient essentiel de notre d\'emonstration est l'obtention de nouvelles estim\'ees multilin\'eaires pour les harmoniques sph\'eriques.
\end{altabstract}

\maketitle
\mainmatter
\section{Introduction}
Let $(M,g)$ be a compact smooth boundary-less Riemannian manifold of dimension $d\geq 2$. Denote by ${\mathbf \Delta}$ the 
Laplace operator associated
to the metric $g$. In the case $d=2$, we discovered in \cite{BGT'03} a bilinear generalization of the well-known Sogge estimates 
\cite{So1,So2,So3} for
$L^p$ ($p\geq 2$) norms of $L^2$ normalized eigenfunctions of ${\mathbf \Delta}$. These bilinear estimates play a central role in the analysis
of \cite{BGT'03} concerning the nonlinear Schr\"odinger equation (NLS) posed on $M$. The goal of this paper is to generalize our
bilinear estimate of \cite{BGT'03} to all higher dimensions and to deduce new results regarding the global existence of solutions for
NLS when $d=3$.

We consider thus the Cauchy problem for NLS
\begin{equation}\label{intr-1}
iu_{t}+{\mathbf \Delta} u=F(u),\quad u|_{t=0}=u_{0}\, .
\end{equation}
In (\ref{intr-1}), $u$ is a complex valued function on $M$. The nonlinear interaction $F$ satisfies $F(0)=0$ and is
supposed of the form 
$
F=\frac{\partial V}{\partial \bar{z}}
$
with $V\in C^{\infty}(\C\,;\,\R)$ satisfying 
\begin{equation}\label{gauge}
V(e^{i\theta}z)=V(z),\quad \theta\in\R, \, z\in\C,
\end{equation}
and, for some $\alpha>1$,
\begin{equation*}
|\partial_{z}^{k_1}\partial_{\bar{z}}^{k_2}\, V(z)|\leq C_{k_1,k_2}(1+|z|)^{1+\alpha-k_1-k_2}\, .
\end{equation*}
The number $\alpha$ involved in the second condition on $V$ corresponds to the ``degree'' of the nonlinearity $F(u)$ in (\ref{intr-1}).
Under these assumptions on $F$, NLS can be seen as a Hamiltonian equation in an infinite dimensional phase space.
It follows from that Hamiltonian structure that smooth solutions of (\ref{intr-1}) enjoy the conservation laws
\begin{equation}\label{intr-2}
\|u(t,\cdot)\|_{L^2}=\|u_0\|_{L^2},\quad E(u(t))=E(u_0),
\end{equation}
where the energy functional $E$ reads as follows,
\begin{equation}\label{intr-3}
E(u)=\int_{M} |\nabla_{g} u|^{2}\,dx +\int_{M}V(u)\,dx \, .
\end{equation}
In view of (\ref{intr-2}) and (\ref{intr-3}), the local well-posedness of (\ref{intr-1}) in 
$H^1(M)$ (with time existence depending upon the $H^1$ norm) is of particular importance.
If for example $V\geq 0$ and $(d-2)\alpha\leq d+2$, (\ref{intr-2}) provides an $H^1$ a priori bound and thus
the local well-posedness of (\ref{intr-1}) in $H^1$ implies the global well-posedness in $H^1$.
Let us notice, on the other hand, 
that the local well-posedness of (\ref{intr-1}) in $H^s$, $s>d/2$ can be obtained by the classical energy 
method (see \cite{Lions}).
If $M$ is two dimensional, the well-posedness of  (\ref{intr-1}) in $H^1(M)$ is established in \cite{BGT1}. In this case, the issue
is to get an improvement of $\varepsilon$ derivatives with respect to the energy method. In \cite{BGT1}, this  $\varepsilon$ gain is 
achieved by a
Strichartz type inequality (with derivative loss). Therefore, for $d=2$, the $H^1$ well-posedness theory for (\ref{intr-1}) is completed.
Moreover, in the recent paper \cite{BGT'03}, we establish a sharp $H^s$ theory in the case $M=S^2$, as far as cubic nonlinearities
are concerned.

In three dimensions, the $H^1$ theory for (\ref{intr-1}) becomes much harder.
In the case $d=3$, 
the Strichartz type inequalities established in
\cite{BGT1} yield the local well-posedness of (\ref{intr-1}) in $H^s$, $s>1$, as far as $\alpha\leq 3$.
Notice that this is already a significant improvement with respect to the energy approach. Unfortunately, it barely misses the crucial
$H^1$ regularity. However, in \cite{BGT1}, we succeeded in using the conservation laws  (\ref{intr-2}) in order to get {\it global}
$H^s$, $s>1$ strong solutions. By ``strong $H^s$ solutions'', we mean the existence, the uniqueness, the propagation of regularity and 
the uniform continuous dependence in bounded subsets of initial data in $H^s$.
Moreover, the methods of \cite{BGT1} also yield uniqueness of $H^1$ weak solutions.

On the other hand, if $M$ is the torus $\T^3$ and $\alpha<5$, the global existence of $H^1$ strong solutions of (\ref{intr-1}) was 
established by Bourgain~\cite{Bo1}.
The approach in \cite{Bo1} is based on an ingenious use of multiple Fourier series and thus relies deeply on the particular
structure of the torus. In this paper, we will prove the counterpart of this result of Bourgain to the cases of the sphere $S^3$
and the product manifold $S^2_{\rho}\times S^1$, where $S^d_{\rho}$, $d\geq 1$ is the 
embedded sphere of radius $\rho$ in $\R^{d+1}$.  
\begin{theoreme}\label{thm1}
Let $M=S^3$ or $M=S^2_{\rho}\times S^1$ endowed with the standard metrics. 
Suppose that $\alpha<5$ and $V(z)\geq -C(1+|z|)^{\beta}$, $\beta<10/3$. Then there exists a space $X$ continuously embedded in 
$C(\R\,;\, H^{1}(M))$ such that for every $u_{0}\in H^1(M)$ there exists a unique global solution $u\in X$ of the Cauchy problem
(\ref{intr-1}). Moreover
\begin{enumerate}
\item
For every $T>0$, the map $u_0\mapsto u\in C([-T,T]\,;\, H^{1}(M))$ is Lipschitz continuous on bounded sets of $H^1(M)$.
\item
If $u_0\in H^{\sigma}(M)$, $\sigma\geq 1$, then for every $t\in\R$, $u(t)\in H^{\sigma}(M)$.
\end{enumerate}
\end{theoreme}
Let us make some comments about this result.
The condition 
$$
V(z)\geq -C(1+|z|)^{\beta},\quad \beta<10/3
$$ 
is classically (see e.g. Cazenave \cite{Caz}) imposed to ensure that the energy controls the $H^1(M)$ norm (defocusing
case). 

The space $X$ will be defined in section~\ref{sec.3} as a local version of Bourgain space $X^{1,b}$. 
It is used to ensure the uniqueness of solutions. However, observe that if  $\sigma> 3/2$, then the uniqueness holds 
in the class $C(\R\,;\, H^{\sigma}(M))$. 
In particular, our theorem implies that for any smooth data $u_0$, there exists a unique global smooth solution.

In the appendix of this paper, we show that Theorem \ref{thm1} can not hold for $\alpha>5$.
The proof is based on an adaptation of an argument of a recent paper of Christ-Colliander-Tao \cite{CCT}
to the setting of compact Riemannian manifolds. The critical case $\alpha =5$ is still open.

Let us recall that the result of Theorem \ref{thm1} was known if we replace $M$ with the Euclidean space $\R^3$ 
(see Ginibre-Velo \cite{GV1} and Kato~\cite{Kato}).
To get the $H^1(\R^3)$ well-posedness of (\ref{intr-1}), for $\alpha<5$, it is sufficient to apply the Picard 
iteration scheme to the Duhamel formulation of f (\ref{intr-1}) in the space
$
L^{2}_{T}W^{1,6}(\R^3)\cap L^{\infty}_{T}H^{1}(\R^3),
$
where $T$ depends only on $\|u_0\|_{H^1}$. The approach on $\R^3$ breaks down in the case of a compact manifold since 
the corresponding Strichartz type estimates have to 
encounter some unavoidable derivative losses (see \cite{Bo1,BGT1,BGT2}). In order to deal with such losses, bilinear improvements
of the Strichartz inequalities are very useful (see e.g. \cite{Bo1,KM,KlMa95,Tao,BGT'03}).
This is the approach that we will adopt in the proof of  Theorem \ref{thm1}
when $M=S^3$. 
The proof in the case $M=S^{2}_{\rho}\times S^{1}$ is more intricate.
The bilinear Strichartz estimates that we are able to prove in the case $M=S^{2}_{\rho}\times S^{1}$ are considerably weaker compared
to the corresponding estimates for $M=S^3$. However, they are sufficient to treat the case $\alpha\leq 4$.
The crucial new point involved in the analysis on $S^{2}_{\rho}\times S^{1}$ is that we can prove a trilinear improvement
of the Strichartz estimate which enables one to treat the case $\alpha=5$ for data in $H^{s}(S^{2}_{\rho}\times S^{1})$, $s>1$.
A suitable interpolation (in the framework of a Littlewood-Paley analysis) between the bilinear and the  trilinear approach
finally completes the argument in the case $M=S^{2}_{\rho}\times S^{1}$.
 
The results of Ginibre-Velo~\cite{GV1} on $\mathbb{R}^3$, of Bourgain~\cite{Bo1} on $\mathbb{T}^3$ 
(and more recently on the irrational three dimensional torus~\cite{Bo4}), and Theorem~\ref{thm1} were obtained for seemingly 
different reasons in each case. For the torus the eigenfunctions enjoy very good algebraic properties and $L^p$ bounds whereas the 
spectrum is ``badly'' localized. On the other hand for the sphere $S^3$, the eigenfunctions present ``bad'' concentration properties 
but the spectrum is very well localized, and the manifold $S^2\times S^1$ has an intermediate behavior. The balance between these 
properties (concentration of eigenfunctions and repartition of the spectrum) leads to the suggestion that a similar result might hold 
for {\em any} three dimensional manifold. The proof of this conjecture would necessitate a general analysis of the 
Schr\"odinger group, unifying these different approaches, which seems to be out of reach at the present moment.

The $H^1$ theory for (\ref{intr-1}) in dimensions $d\geq 4$ remains an open problem. The only known result in this direction is that
of Bourgain \cite{Bo2} who gets {\it global} $H^s(\T^4)$ solutions, if $\alpha\leq 2$, $s>1$.

It seems that the obstructions to extending our approach to high dimensions are not only of technical nature since in \cite{BGT2}
we have shown that for no $\alpha\in]1,2]$ (even very close to $1$), the Cauchy problem (\ref{intr-1}), posed on $S^6$ can
have strong $H^1$ solutions in the sense explained above. Interestingly, the result of \cite{BGT2} is in strong contrast with
the situation on $\R^6$ (see \cite{BGT-lund}).
\par
We now turn to the crucial step in the proof of Theorem~\ref{thm1}. To that purpose, we introduce the following notation : given
$\nu\geq 1$, we set
\begin{equation*}
\Lambda(d,\nu): = 
\begin{cases}
\nu^{\frac{1}{4}} &\text{ if $d=2$}\\
\nu^{\frac{1}{2}}\log^{1/2} (\nu) &\text{ if $d=3$}\\
\nu^{\frac{d-2}{2}} &\text{ if $d\geq 4$}.
\end{cases}
\end{equation*}
With this notation, we have the following multilinear eigenfunction estimates.
\begin{theoreme}\label{thm2}
There exists $C>0$ such that, if $H_p$ and $H_q$ are two spherical harmonics of respective degrees $p$ and $q$,  
\begin{equation}\label{thm2-1}
\|H_p H_q\|_{L^{2}(S^d)}
\leq C \Lambda(d, \min (p,q)+1) \|H_p\|_{L^{2}(S^d)}\|H_q\|_{L^{2}(S^d)}\, .
\end{equation}
Moreover for any  $p\geq q\geq r \geq 0$, the following trilinear estimates hold
\begin{equation}\label{thm2-2}
\|H_p H_q H_r\|_{L^{2}(S^2)}
\leq
C[(1+q)(1+r)]^{\frac{1}{4}}
\|H_p\|_{L^{2}(S^2)}\|H_q\|_{L^{2}(S^2)}\|H_r\|_{L^{2}(S^2)}\, .
\end{equation}
Estimates (\ref{thm2-1}) and (\ref{thm2-2}) are sharp, apart from the logarithmic loss in (\ref{thm2-1}) for $d=3$.
\end{theoreme}
\begin{remarque}
As an easy consequence of (\ref{thm2-1}), one can prove the corresponding estimate to (\ref{thm2-2}), for $d\geq 3$,
\begin{equation}\label{thm2-2-novo}
\|H_p H_q H_r\|_{L^{2}(S^d)}
\leq
C \Lambda(d, q+1)(1+r)^{\frac{d-1}{2}}
\|H_p\|_{L^{2}(S^d)}\|H_q\|_{L^{2}(S^d)}\|H_r\|_{L^{2}(S^d)}\, .
\end{equation}
Indeed it suffices to use that
the $L^{\infty}(S^d)$ norm of $H_r$ is bounded by $(1+r)^{\frac{d-1}{2}}$ (Weyl bound) and  (\ref{thm2-1}) for 
the product $H_p H_q$. 

In view of further possible developments, we will also prove in section~\ref{sec.2} that for every $\eta\in ]0,1]$ there exists 
$C_{\eta}$ such that
\begin{equation}\label{thm2-2-novo-bis}
\|H_p H_q H_r\|_{L^{2}(S^3)}
\leq C_\eta (1+q)^{\frac{1}{2}+\eta}
(1+r)^{1-\eta}
\|H_p\|_{L^{2}(S^3)}\|H_q\|_{L^{2}(S^3)}\|H_r\|_{L^{2}(S^3)}\, .
\end{equation}
\end{remarque}

In fact, we deduce Theorem~\ref{thm2} as a consequence of a more general statement concerning the approximated spectral projectors 
$\chi(\sqrt{-{\mathbf \Delta}}-\lambda)$, $\lambda\gg 1$, $\chi\in {\cal S}(\R)$, where ${\mathbf \Delta}$ is the Laplace operator on
an arbitrary compact Riemannian manifold $(M,g)$ (see Theorem \ref{thm3} below).

Notice that when $p=q=r$, apart from the $\log$ loss in $3d$, we recover some particular case of the $L^p-L^2$ {\it linear}
estimates of Sogge \cite{So1,So2,So3}. In the proof of Theorem~\ref{thm1}, we typically apply Theorem \ref{thm2} 
for $p\gg q$ and thus estimates (\ref{thm2-1}),  (\ref{thm2-2}) are used in their full strength.

In the case $d=2$, estimate (\ref{thm2-1}) has already appeared in our previous paper \cite{BGT'03}.
In \cite{BGT'03}, the proof is inspired by H\"ormander's work \cite{Ho1} on Carleson-Sj\"olin type operators.
The proof we present here is different even for $d=2$ and relies on a ``bilinearization'' of the arguments in 
\cite{So1,So2,So3}. After several preliminaries, we reduce the matters to two micro-local {\it linear} estimates of quite a
different nature. The first one is applied to the higher frequency eigenfunction and is in the spirit of the $L^2$ 
boundedness of spectral projectors. The second one is applied to the smaller frequency eigenfunctions and relies on a dispersive
(curvature) effect. As far as the optimality of (\ref{thm2-1}), (\ref{thm2-2}) is concerned, we notice that it is achieved either 
by testing the estimates against eigenfunctions concentrating on an equator or by testing against zonal eigenfunctions
concentrating on a point.\par
Let us mention that estimates (\ref{thm2-1}), (\ref{thm2-2}) and a sketch of the proof of (\ref{thm2-1}) appeared in \cite{BGT-note}.
\par
The rest of this paper is organized as follows. In section~\ref{sec.2} we prove Theorem~\ref{thm2}. 
In section~\ref{sec.3} we set up the framework of Bourgain's spaces and reduce the proof of Theorem~\ref{thm1} 
to the obtaining of nonlinear estimates in this framework. Section~\ref{sec.4} consists in two parts. 
First we prove bilinear Strichartz estimates for the linear Schr\"odinger group on $S^3$. 
Then we show that Theorem~\ref{thm1} holds for any three dimensional manifold on which these estimates are true. 
Section~\ref{sec.5} also consists in two parts. First we prove trilinear Strichartz estimates for the linear Schr\"odinger 
group on the product manifold $S^2_\rho\times S^1$ and then we show that Theorem~\ref{thm1} holds for any three dimensional manifold 
on which these estimates are true. An appendix is devoted to the proof of the optimality of the quintic threshold.
\par
{\bf Acknowledgements.} 
We are grateful to J. Bourgain for sending us his manu\-script~\cite{Bo4} and 
H. Koch and D. Tataru for interesting discussions about spectral projectors.
\section{Multilinear eigenfunction estimates}
\label{sec.2}
In this section we prove Theorem \ref{thm2}, and more generally the corresponding result for spectral projectors on arbitrary
compact manifolds.
\subsection{On the optimality of the estimates} 
We first consider the optimality of (\ref{thm2-1}) in the case $d=2,3$. Let us see $S^d$ as a hyper-surface in $\R^{d+1}$, i.e.
$$
S^{d}=\{(x_1,\dots,x_{d+1})\in\R^{d+1}\,:\, x_{1}^{2}+\cdots +x_{d+1}^{2}=1\}\, .
$$
Let us define the highest weight spherical harmonics $R_{p}=(x_1+ix_2)^{p}$ which concentrate, for $p\gg 1$, on the closed geodesic
(a big circle) $x_1^2+x_2^2=1$. An easy computation shows that
$$
\|R_{p}\|_{L^{2}(S^d)}\approx p^{-\frac{d-1}{4}},\quad p\gg 1.
$$
Clearly $R_{p}R_{q}=R_{p+q}$ and therefore there exist constants $C$, $\widetilde{C}$ such that for every $(p,q)$,
$$
\|R_{p}R_{q}\|_{L^2(S^d)}\geq C (p+q)^{-\frac{d-1}{4}}
\geq \widetilde{C}(\min(p,q))^{\frac{d-1}{4}}\|R_{p}\|_{L^2(S^d)}\|R_{q}\|_{L^2(S^d)}\, .
$$
Therefore, for $d=2,3$, estimate (\ref{thm2-1}) turns out to be optimal, modulo the logarithmic loss in $3d$.
In the same way, since $R_{p}R_{q}R_{r}=R_{p+q+r}$, estimate (\ref{thm2-2}) is optimal by testing it on $R_{p}$, $R_{q}$ and $R_{r}$. 

Let us now consider the case $d\geq 4$. In this case the optimality of (\ref{thm2-1}) is given by the zonal spherical harmonics.
Let us a fix a pole on $S^d$. If we consider functions on $S^d$ depending only on the geodesic
distance to the fixed pole, we obtain the zonal functions on  $S^d$. The zonal functions can be expressed in terms of zonal
spherical harmonics which in their turn can be expressed in terms of the classical Jacobi polynomials (see e.g. \cite{So1}).
Using asymptotics for the  Jacobi polynomials (see \cite{Sz},\cite[Lemma 2.1]{So1}) we can obtain the following 
representation for the zonal
spherical harmonics $Z_{p}$ of degree $p$, in the coordinate $\theta$, 
\begin{equation}\label{szego}
Z_{p}(\theta)
=
C(\sin\theta)^{-\frac{d-1}{2}}
\Big\{
\cos[(p+\alpha)\theta+\beta]+\frac{{\cal O}(1)}{p\, \sin\theta}
\Big\},\quad \frac{c}{p}\leq \theta\leq\pi-\frac{c}{p}\, ,
\end{equation}
where $\alpha$ and $\beta$ are some fixed constants depending only on $d$.
Moreover, we have a point-wise concentration 
\begin{equation}\label{point}
|Z_{p}(\theta)|\approx  p^{\frac{d-1}{2}},\quad \theta\notin [c/p,\pi-c/p]
\end{equation}
and $\|Z_{p}\|_{L^2(S^d)}\approx 1$. Let $q\gg p$. Then 
$$
\|Z_{p} Z_{q}\|_{L^2(S^d)}^{2}=\int_{0}^{\pi}Z_{p}^{2}(\theta) Z_{q}^{2}(\theta) (\sin\theta)^{d-1} d\theta
\geq
\int_{c/q}^{c/p}Z_{p}^{2}(\theta) Z_{q}^{2}(\theta) (\sin\theta)^{d-1} d\theta\, .
$$
Using (\ref{point}), we get
$$
\|Z_{p} Z_{q}\|_{L^2(S^d)}^{2}
\geq Cp^{d-1}\int_{c/q}^{c/p}Z_{q}^{2}(\theta) (\sin\theta)^{d-1} d\theta.
$$
In view of (\ref{szego}),
$$
\|Z_{p} Z_{q}\|_{L^2(S^d)}^{2}
\geq
Cp^{d-1}[I_1-I_2],
$$
where
$$
I_{1}=\int_{c/q}^{c/p}\cos^{2}[(p+\alpha)\theta+\beta]d \theta\geq \frac{C}{p}
\quad
{\rm and }
\quad
I_{2}=\frac{1}{q^2}\int_{c/q}^{c/p}\frac{1}{(\sin\theta)^2}d\theta\leq \frac{C}{q}\ll \frac{C}{p}\, .
$$
Therefore
\begin{equation}\label{lll1}
\|Z_{p} Z_{q}\|_{L^2(S^d)}^{2}\geq Cp^{d-2}
\|Z_{p}\|_{L^2(S^d)}^{2}\|Z_{q}\|_{L^2(S^d)}^{2},
\end{equation}
if  $p\ll q$. 
Let finally $p\approx q$. Using (\ref{point}), we get
\begin{equation}\label{lll2}
\|Z_{p} Z_{q}\|_{L^{2}}^{2}
\geq  
C\,p^{2(d-1)}
\int_{0}^{c/p}(\sin\theta)^{d-1}d\theta
\geq 
\tilde{C}\,p^{2(d-1)}\, p^{-d}
=\tilde{C}\,p^{d-2}\, .
\end{equation}
Therefore, collecting (\ref{lll1}) and (\ref{lll2}), we obtain
$$
\|Z_{p} Z_{q}\|_{L^2(S^d)}\geq C(\min(p,q))^{\frac{d-2}{2}}
\|Z_{p}\|_{L^2(S^d)}\|Z_{q}\|_{L^2(S^d)}
$$
which proves the optimality of (\ref{thm2-1}), for $d\geq 3$, modulo the logarithmic loss in $3d$.
Let us finally notice that similarly we can prove that for $p\geq q\geq r$
$$
\|Z_{p} Z_{q} Z_{r}\|_{L^2(S^d)}\geq Cq^{\frac{d-2}{2}}r^{\frac{d-1}{2}}
\|Z_{p}\|_{L^2(S^d)}\|Z_{q}\|_{L^2(S^d)}\|Z_{r}\|_{L^2(S^d)}
$$
which proves the optimality of (\ref{thm2-2-novo}), for $d\geq 3$, apart from the logarithmic loss in $3d$, and
the optimality of (\ref{thm2-2-novo-bis}) apart from the $\eta$ shift.
\subsection{A first reduction}
Let $(M,g)$ be a compact smooth Riemannian manifold without boundary of dimension $d$  and 
${\mathbf \Delta}$ be the Laplace operator on functions on $M$. 
It turns out that estimates (\ref{thm2-1}),  (\ref{thm2-2}) and  (\ref{thm2-2-novo-bis}) 
can be deduced from the following more general result.
\begin{theoreme}\label{thm3}
Let $\chi\in {\cal S}(\R)$. For $\lambda \in \R$,  denote by  
$\chi_\lambda = \chi(\sqrt{-{\mathbf \Delta}}-\lambda)$ the approximated spectral projector around $\lambda$.
 There exists $C$ such that for any $\lambda, \mu \geq 1 $, $f,g \in L^{2}(M)$,
\begin{equation}\label{eqbilin}
\|\chi_\lambda f \,\chi_\mu g
\|_{L^{2}(M)}
\leq C \Lambda(d, \min (\lambda, \mu)) \|f\|_{L^{2}(M)}\|g\|_{L^{2}(M)}\, .
\end{equation}
Moreover, in the case $d=2$, for any  $1\leq \lambda\leq \mu\leq \nu $, $f,g, h \in L^{2}(M)$, the following trilinear estimate holds
\begin{equation}\label{eqtrilin}
\|\chi_{\lambda}f\,  \chi_{\mu}g \, \chi_{\nu}h
\|_{L^{2}(M)}
\leq
C(\lambda\mu)^{\frac{1}{4}}\|f\|_{L^{2}(M)}\|g\|_{L^{2}(M)}\|h\|_{L^{2}(M)}.
\end{equation}
Finally, in the case $d=3$, for any  $1\leq \lambda\leq \mu\leq \nu$, $f,g, h \in L^{2}(M)$, $\eta\in ]0,1]$,
the following trilinear estimate holds
\begin{equation}\label{eqtrilin-3}
\|\chi_{\lambda}f\,  \chi_{\mu}g \, \chi_{\nu}h\|_{L^{2}(M)}\leq
C_{\eta}\lambda^{1-\eta}\mu^{\frac{1}{2}+\eta}\|f\|_{L^{2}(M)}\|g\|_{L^{2}(M)}\|h\|_{L^{2}(M)}.
\end{equation}
\end{theoreme}
\begin{remarque}
If one is only interested in estimates for single eigenfunctions,
the bounds provided by Theorem \ref{thm3} seem to be relevant for ``sphere like manifolds''
but they are far from the optimal ones in the case of the torus. For example, the classical result of
Zygmund \cite{Z} says that there exists a constant $C$ such that for every couple $(f,g)$ of eigenfunctions of the Laplace
operator on the torus $\T^2$, one has 
$$
\|f\,g\|_{L^2(\T^2)}\leq C\|f\|_{L^2(\T^2)}\|g\|_{L^2(\T^2)}\, .
$$
We refer to Bourgain \cite{Bo3} for further extensions of Zygmund's result.
\end{remarque}
A first reduction in the proof of Theorem \ref{thm3} is that it suffices to prove it for one fixed non trivial function $\chi$.
\begin{lemme}\label{chi}
Suppose that the assertion of Theorem \ref{thm3} holds for a bump function $\chi\in{\cal S}(\R)$ which is not identically zero.
Then it holds for any other choice of the bump function.
\end{lemme}
\begin{proof}
Suppose that  (\ref{eqbilin}) holds for a nontrivial $\chi\in{\cal S}(\R)$.
Then, there exists $x_0\in \R$ such that $\chi(x_0)\neq 0$ and moreover there exists $\delta>0$ such that $\chi(x)\neq 0$ for $x\in\R$
satisfying $|x-x_0|<2\delta$. Using a partition of unity argument, we can find $\psi\in C_{0}^{\infty}(\R)$
supported in 
$
\{x\in\R\,:\, |x|<\frac{3\delta}{4}\}
$
such that
\begin{equation}\label{pu}
\sum_{n\in \Z}\psi(x-n\delta)=1\, .
\end{equation}
Thanks to the support properties of $\psi$ and $\chi$, we can write
\begin{equation}\label{div}
\psi(x-n\delta-\lambda)=\chi(x+x_0-n\delta-\lambda)\, \frac{\psi(x-n\delta-\lambda)}{\chi(x+x_0-n\delta-\lambda)}\, .
\end{equation}
Notice that the second factor in the right hand-side of (\ref{div}) is uniformly bounded. Therefore, using that  (\ref{eqbilin})
holds for $\chi$, we obtain the estimate
\begin{multline}\label{loc-psi}
\|\psi(\sqrt{-{\mathbf \Delta}}-n\delta-\lambda)(f)\, \psi(\sqrt{-{\mathbf \Delta}}-m\delta-\mu)(g)\|_{L^2}
\leq
\\
\leq
C\Lambda(d, \min (|n|+\lambda, |m|+\mu))
\|f\|_{L^{2}}\|g\|_{L^{2}}\, .
\end{multline}
Let us now take an arbitrary function $\chi_{1}\in {\cal S}(\R)$. Using (\ref{pu}), we can write
\begin{equation}\label{eee}
\chi_{1}(\sqrt{-{\mathbf \Delta}}-\lambda)f
=\sum_{n\in\Z}
\psi(\sqrt{-{\mathbf \Delta}}-n\delta-\lambda)\,\chi_{1}(\sqrt{-{\mathbf \Delta}}-\lambda)f \, .
\end{equation}
Let $\tilde{\psi}\in C_{0}^{\infty}(\R)$ be equal to one on the support of $\psi$. Then clearly
\begin{equation}\label{ddd}
|\chi_{1}(x-\lambda)\tilde{\psi}(x-\lambda-n\delta)|
\leq
\frac{C_N}{(1+|x-\lambda|)^{N}(1+|x-\lambda-n\delta|)^{N}}
\leq
\frac{\widetilde{C}_N}{(1+|n|)^{N}}\, .
\end{equation}
Using the expansion (\ref{eee}) together with (\ref{loc-psi}) and (\ref{ddd}) yields
\begin{multline*}
\|\chi_{1}(\sqrt{-{\mathbf \Delta}}-\lambda)(f)\, \chi_{1}(\sqrt{-{\mathbf \Delta}}-\mu)(g)\|_{L^2}
\\
\leq 
\sum_{(n,m)\in \Z^2}
\frac{C_N \Lambda(d, \min (|n|+\lambda, |m|+\mu))}{(1+|n|)^{N}(1+|m|)^{N}}\|f\|_{L^{2}}\|g\|_{L^{2}}
\\ 
\leq
C\Lambda(d, \min (\lambda, \mu))\|f\|_{L^{2}}\|g\|_{L^{2}}\, .
\end{multline*}
Hence  (\ref{eqbilin}) holds for $\chi_1$. The proof of the independence of (\ref{eqtrilin}) 
and (\ref{eqtrilin-3}) with respect to the bump function $\chi$ is very similar and thus we will omit it.
\end{proof}
\subsection{Reduction to oscillatory integral estimates and main properties of the phase function}
Following \cite[Chap. 4]{So3}, thanks to Lemma \ref{chi}, it is sufficient to prove Theorem \ref{thm3} with $\chi$ such
that $\widehat{\chi}(\tau)$ is supported in the set
$$
\{\tau\in\R\,:\, \varepsilon\leq \tau\leq 2\varepsilon\},
$$ 
where $\varepsilon>0$ is a {\it small}
number to be determined later. We can write 
$$
\chi_{\lambda}f=\frac{1}{2\pi}\int_{\varepsilon}^{2\varepsilon}
e^{-i\lambda\tau}\widehat{\chi}(\tau)(e^{i\tau\sqrt{-{\mathbf \Delta}}}f)d\tau\, .
$$
For $\varepsilon\ll 1$ and $|\tau|\leq 2\varepsilon$, 
using a partition of the unity on $M$, we can represent $e^{i\tau\sqrt{-{\mathbf \Delta}}}$ as a Fourier integral operator
(see e.g. \cite{Ho-1}). Therefore $\chi_{\lambda}$ can also be represented as such. After a stationary phase
argument (see \cite[Chap. 5]{So3}) we can represent $\chi_{\lambda}f$ as follows.
\begin{lemme}\label{sogge}
There exists $\varepsilon_{0}>0$ such that for every $\varepsilon\in]0,\varepsilon_{0}[$, every $N\geq 1$, we have the splitting
\begin{equation}\label{reste}
\chi_\lambda f=\lambda^{\frac {d-1} 2} T_{\lambda}f+R_{\lambda}f, 
\end{equation}
with
$$
\|R_{\lambda}f\|_{H^{k}(M)}\leq C_{N,k}\lambda^{k-N}\|f\|_{L^2(M)},\quad \, k=0,\dots,N\, .
$$
Moreover there exist $\delta>0$ and, for every $x_0\in M$, a system of coordinates $V\subset\R^d$, 
containing $0\in \R^d$ such for $x\in V$, $|x|\leq \delta$,
$$
T_{\lambda}f(x)=\int_{\R^d}e^{i\lambda\varphi(x,y)}a(x,y, \lambda)f(y)dy
$$
where $a(x,y,\lambda)$ is a polynomial in $\lambda^{-1}$ with smooth coefficients supported in the set 
\begin{equation*}
\{(x,y)\in V\times V\, : \, |x| \leq \delta \ll \frac{\varepsilon}{C}\leq |y|\leq C\varepsilon\}
\end{equation*}
and   $-\varphi(x,y)=d_{g}(x,y)$ is the geodesic distance between $x$ and $y$.
\end{lemme}
\begin{remarque}
Let us notice that one can use $\chi(-\lambda^{-1}{\mathbf \Delta}-\lambda)$ as approximated spectral projector 
instead of $\chi(\sqrt{-{\mathbf \Delta}}-\lambda)$.
In that case one should use semi-classical calculus for the approximation of $\exp(it\lambda^{-1}{\mathbf \Delta})$, $\lambda\gg 1$,
as we did in \cite{BGT1}. 
\end{remarque}
In view of Lemma \ref{sogge}, to prove  (\ref{eqbilin}), it is enough to show
\begin{equation}\label{eqbilin-pak}
\|{T}_{\lambda}f\,{T}_{\mu}g\|_{L^2}
\leq
C\Lambda(d,\lambda)(\lambda \mu)^{-\frac{d-1}{2}}
\|f\|_{L^2}\|g\|_{L^2},
\end{equation}
uniformly for $1\leq\lambda\leq \mu$. Indeed, using (\ref{reste}), one has to evaluate in $L^2$ the products 
$$
T_{\lambda}f\,R_{\mu}g,\quad R_{\lambda}f\,T_{\mu}g\quad R_{\lambda}f\,R_{\mu}g.
$$
The products involving $R_{\mu}$ are straightforward to estimate while for $R_{\lambda}f\,T_{\mu}g$,
using the $L^2$ boundedness of $\chi_{\mu}$, we write
$$
\|R_{\lambda}f\,T_{\mu}g\|_{L^2}\leq C \|R_{\lambda}f\|_{L^{\infty}}\|T_{\mu}g\|_{L^2}
\leq C_{N}\lambda^{-N}
\mu^{-\frac{d-1}{2}}\|f\|_{L^2}\|g\|_{L^2}\, .
$$
Furthermore, we notice that once~\eqref{eqbilin-pak} is proved (at least for $d=2$), to prove (\ref{eqtrilin}) it is enough to show that for $d=2$,
\begin{equation*}
\|{T}_{\lambda}f\,{T}_{\mu}g\, {T}_{\nu} h \|_{L^2}
\leq
C(\lambda\mu)^{-\frac{1}{4}}\nu^{-\frac{1}{2}}\|f\|_{L^2}\|g\|_{L^2}\|h\|_{L^2},
\end{equation*}
uniformly for $1\leq\lambda\leq \mu\leq \nu$. In this case there are more remainder terms to estimate.
The most difficult one is ${R}_{\lambda}f\,{T}_{\mu}g\, {T}_{\nu} h $. This term can be evaluated, by using
(\ref{eqbilin-pak}) for $d=2$, as follows
$$
\|{R}_{\lambda}f\,{T}_{\mu}g\, {T}_{\nu} h\|_{L^2}
\leq
\|{R}_{\lambda}f\|_{L^{\infty}}
\|{T}_{\mu}g\, {T}_{\nu} h\|_{L^2}
\leq C_{N}\lambda^{-N}\mu^{-\frac{1}{4}}\nu^{-\frac{1}{2}}
\|f\|_{L^2}\|g\|_{L^2}\|h\|_{L^2}\, .
$$
Similarly, to prove
(\ref{eqtrilin-3}), it is enough to show that for $d=3$,
\begin{equation*}
\|{T}_{\lambda}f\,{T}_{\mu}\,g {T}_{\nu}\, h \|_{L^2}
\leq
C\lambda^{-\eta}\mu^{-\frac{1}{2}+\eta}\nu^{-1}
\|f\|_{L^2}\|g\|_{L^2}\|h\|_{L^2},
\end{equation*}
uniformly for $1\leq\lambda\leq \mu\leq \nu$.
In this case, we estimate ${R}_{\lambda}f\,{T}_{\mu}g\, {T}_{\nu} h $, 
by another use of (\ref{eqbilin-pak}), as follows
\begin{multline*}
\|{R}_{\lambda}f\,{T}_{\mu}g\, {T}_{\nu} h\|_{L^2}
\leq
\|{R}_{\lambda}f\|_{L^{\infty}}
\|{T}_{\mu}g\, {T}_{\nu} h\|_{L^2}
\leq
\\
\leq
C_{N}\lambda^{-N}
\log^{1/2}(\mu)\mu^{-\frac{1}{2}}\nu^{-1}
\|f\|_{L^2}\|g\|_{L^2}\|h\|_{L^2}
\leq 
\\
\leq
C_{N,\eta}\lambda^{-N}\mu^{-\frac{1}{2}+\eta}\nu^{-1}
\|f\|_{L^2}\|g\|_{L^2}\|h\|_{L^2}\, ,
\end{multline*}
where $\eta>0$.
\vskip .5cm
Next, we represent $y$ in geodesic (polar) coordinates as $y=\exp_{0}(r\omega)$, $r >0$, $\omega \in S^{d-1}$.
For $|x|\leq \delta$ and $\omega\in S^{d-1}$, we define the frozen phase $\varphi_{r}$,
$$
\varphi_r(x, \omega)=\varphi(x,\,\exp_{0}(r\omega)) \, .
$$
We now state the main property of the phase $\varphi_r$.
\begin{lemme}\label{pak}
There exists $\varepsilon>0$ such that
for every $r\in [\varepsilon/C, C\varepsilon]$, every 
$$
\omega=(\omega_1,\dots,\omega_{d})\in S^{d-1}\subset \R^{d},
$$
we have the identity,
$$
\nabla_{x}\varphi_{r}(0,\omega)=\omega\, .
$$
\end{lemme}
\begin{proof}
The proof for $d=2$ is given in \cite{BGT'03}. The extension to an arbitrary $d$ is straightforward as we explain below.
For $\varepsilon\ll 1$,
let $y= \exp_0(r \omega)$, $r= -\varphi(0,y)$ and $u=u(x,y)\in {T}_{y}M$ be the unique unit vector in the tangent space to $M$ at 
$y$ such that
$$
\exp_{y}(-\varphi(x,y)u(x,y))=x.
$$
Differentiating with respect to $x$ this identity, we get for $x=0$, and any $h\in T_{0}M$,
\begin{multline}\label{geom}
h=
-g_{0}\Big(\nabla_{x}\varphi(0,y)\,,\,h\Big)\,T_{ru(0,y)}(\exp_{y})\cdot u(0,y)
\\
+T_{ru(0,y)}(\exp_{y})\big(r\, T_{x}u(0,y)\cdot h\big),
\end{multline}
where $T$ denotes the tangential map.

On the other hand, we have
\begin{equation}\label{jacoby}
T_{ru(0,y)}(\exp_{y})\cdot u(0,y)=-\omega\quad \text{ or }\quad u(0,y)= -T_{r\omega} (\exp_0) (\omega).
\end{equation}
Consequently, using Gauss' Lemma (see \cite[3.70]{GHL}), we get
\begin{equation}\label{geom3}
 g_{0}\Big(T_{ru(0,y)}(\exp_{y})\big(r\, T_{x}u(0,y)\cdot h\big)\,,\, \omega\Big)=0
\, .
\end{equation}
Let us now take the scalar product of (\ref{geom}) with $\omega$. Collecting (\ref{geom}), (\ref{geom3}) and~\eqref{jacoby}  yield
$$
g_{0}\big(\omega\,,\,h\big)=g_{0}\big(\nabla_{x}\varphi(0,y)\,,\,h\big),\quad\forall h\in T_{0}M
$$
which completes the proof of Lemma \ref{pak}.
\end{proof} 
Let us notice that there exists a smooth positive function $\kappa(r,\omega)$ such that $dy = \kappa(r,\omega) dr d\omega$.
For $r\in  [\frac{\varepsilon}{C}, C\varepsilon]$ and $\lambda\geq 1$,
we define the operator $T_{\lambda}^{r}$, acting on functions on $S^{d-1}$ via the identity
$$
(T_{\lambda}^{r}f)(x)=\int_{S^{d-1}}
e^{i\lambda\varphi_{r}(x,\omega)}
a_r(x,\omega,\lambda)f(\omega)d\omega,
$$
where $a_r(x, \omega,\lambda) = \kappa(r,\omega)a(x,\text{exp}_0(r\omega),\lambda)$.
Then clearly
$$
(T_{\lambda}f)(x)=\int_{0}^{\infty}(T_{\lambda}^{r}f_r)(x) dr,
$$
where $f_r(\omega) = f(r,\omega)$. Similarly, with $g_q(\omega)=g(q,\omega)$, 
$$
(T_{\lambda}f\,{T}_{\mu}g)(x)
=
\int_{\varepsilon/C}^{C\varepsilon}
\int_{\varepsilon/C}^{C\varepsilon}
(T_{\lambda}^{r}f_r)(x)\,
({T}_{\mu}^{q}g_q)(x)drdq,
$$
and the Minkowski inequality shows that (\ref{eqbilin}) will be a consequence of
\begin{equation}\label{edno}
\|{T}_{\lambda}^{r}f\,{T}_{\mu}^{q}g\|_{L^2}
\leq
C\Lambda(d,\lambda)(\lambda \mu)^{-\frac{d-1}{2}}
\|f\|_{L^2(S^{d-1})}\|g\|_{L^2(S^{d-1})},
\end{equation}
uniformly for $1\leq\lambda\leq \mu$ and $r,q\in [\frac{\varepsilon}{C}, C\varepsilon]$.

Similarly, to prove (\ref{eqtrilin}), it is enough to show
\begin{equation}\label{dve}
\|{T}_{\lambda}^{r}f\,{T}_{\mu}^{q}g\, {T}_{\nu}^{s} h \|_{L^2}
\leq
C(\lambda\mu)^{-\frac{1}{4}}\nu^{-\frac{1}{2}}
\|f\|_{L^2(S^{1})}\|g\|_{L^2(S^{1})}\|h\|_{L^2(S^{1})},
\end{equation}
uniformly for $1\leq\lambda\leq \mu\leq \nu$ and $r,q,s \in [\frac{\varepsilon}{C}, C\varepsilon]$.

Finally, to prove (\ref{eqtrilin-3}), it is enough to show
\begin{equation}\label{tri}
\|{T}_{\lambda}^{r}f\,{T}_{\mu}^{q}g \, {T}_{\nu}^{s} h \|_{L^2}
\leq
C
\lambda^{-\eta}\mu^{-\frac{1}{2}+\eta}\nu^{-1}
\|f\|_{L^2(S^{2})}\|g\|_{L^2(S^{2})}\|h\|_{L^2(S^{2})},
\end{equation}
uniformly for $1\leq\lambda\leq \mu\leq \nu$ and $r,q,s \in [\frac{\varepsilon}{C}, C\varepsilon]$.

Fix a point $\underline{\omega} \in S^{d-1}$. The set 
$$
S_{x}= \{ \nabla_{x} \varphi_r (x,\omega),\quad  \omega\in S^{d-1},\,\, \omega\sim \underline{\omega}\}
$$
is a smooth hyper-surface in $\R^d$. Indeed assuming for instance $\underline{\omega}=(1, 0, \dots, 0)$, then 
$(w_1= \omega _2, \dots, , w_{d-1}= \omega_d)$ is a system of coordinates on $S^{d-1}$ and according to 
Lemma~\ref{pak}, $\nabla_w \nabla_x \varphi_r $ has rank $d-1$.

Following Stein \cite{Stein} and Sogge \cite{So3}, we now state the crucial curvature property.
\begin{lemme}\label{lem.courb}
The hyper-surface $S_{x}$ has non-vanishing principal curvatures:
for $w\in \R^{d-1}$ a local coordinate system near  $\underline{\omega}\in S^{d-1}$, 
if we denote by $\pm n (x,w)$ the normal unit vectors to the surface $S_{x}$ at the point $\nabla_{x}\varphi_r (x,w)$,
then for $x$ close to $0$,
\begin{equation}\label{eq1.2}
\left|\det_{i,j}\big\langle
\frac{\partial^2 }{ \partial w_{j}\partial{w_{i}} } 
\nabla_{x}\varphi_{r}(x,w), n( x,w)\big\rangle\right|\geq c> 0 .
\end{equation}
\end{lemme}
\begin{proof}
The relation (\ref{eq1.2}) is equivalent to the fact that 
$$
w\mapsto n(x,w) \in S^{d-1}
$$ is a local diffeomorphism. Indeed, dropping the $x$ variable for conciseness and denoting by
$$M(w)= \nabla_x \varphi_r(x,w), \qquad n(w) = n(x,w), $$
we have
$$\langle \frac{\partial M} { \partial w_i}, n(w)\rangle =0 \Rightarrow \langle \frac {\partial^2 M}{ \partial w_i \partial w_j}, n \rangle = - \langle \frac{ \partial M } {\partial w_i}, \frac {\partial n} {\partial w_j}\rangle$$
As a consequence, the determinant in~\eqref{eq1.2} is non vanishing if and only if the system of vectors $\frac { \partial n} { \partial w_j}$ is of maximal rank in $T_w S_x$. We deduce that~\eqref{eq1.2} is independent of the choice of coordinates $w$ and it suffices 
to prove it for a particular choice of a coordinate system near $\underline{\omega}$. We can suppose that 
$\underline{\omega}= (1,0, \dots, 0)$ and we choose as coordinates 
$$
w=(w_1,\dots,w_{d-1}):=(\omega_2,\dots,\omega_{d})
$$
We can also assume that at the point $(x=0)$, the metric is diagonal, $g_{i,j}= \delta_{i,j}$.
Using  Lemma \ref{pak}, we get
\begin{equation}\label{eq1.2bis}
\big\langle
\frac{\partial^2 }{ \partial w_{j}\partial{w_{i}} } 
\nabla_{x}\varphi_{r}(0,w), n( 0,0)\big\rangle\mid_{w=0}  =\text{Id}
\end{equation}
and consequently~\eqref{eq1.2} follows by continuity.
\end{proof}
 Denote by $(T^\nu_r)^*$ the formal adjoint of $T^\nu_r$. The kernel of the operator $T^\nu_r ( T^\nu_r)^*$, $K(x,x')$, 
is given by the relation
$$
K(x,x')= \int e^{i\nu ( \varphi_r(x, w) - \varphi_r(x', w))} a_r(x,w,\nu) \overline{a} _r(x', w, \nu)dw\, .
$$
The curvature property of the phase $\varphi_r$ in Lemma~\ref{lem.courb} implies a dispersion inequality for the kernel $K$.
\begin{lemme}\label{lem.dispersion}
There exist $C>0$ such that for any $\nu \geq 1$, 
\begin{equation}\label{eq.dispersion}
|K(x,x')| \leq \frac C { (1+ \nu |x-x'|)^{\frac {d-1} 2 }}\, .
\end{equation}
\end{lemme}
\begin{proof}
Let us write a Taylor expansion
$$
\varphi_{r} (x, w) - \varphi_{r}(x', w) = \langle x-x'\, ,\, \psi(x,x',w) \rangle 
$$
where 
$$
\psi(x,x',w)= \int_0^1 \nabla_x\varphi_{r} (x'+ \theta (x-x'), w) d \theta.
$$
With $\sigma = \frac {x-x' } { |x-x'|}$, we can write 
$$
\varphi_r(x, w) - \varphi_r (x', w) = |x-x'| \Phi (x,x', \sigma, w) 
$$
where 
$$
\Phi (x,x', \sigma, w)= \langle \sigma \, , \, \psi(x,x',w)\rangle
$$
Now we want to prove, with $\lambda = \nu |x-x'|$,
$$ 
|\widetilde{K} (x, x', \sigma) | \leq \frac C { (1+ \lambda)^{\frac {d-1 } 2 }}
$$
where 
\begin{equation}\label{eq.kernel}
\widetilde{K}(x,x', \sigma)= \int e^{i\lambda \Phi(x,x', \sigma, w)} a_r(x,w,\nu) \overline{a} _r(x', w, \nu)dw\, .
\end{equation}
From the definition of the normal $n(x, w)$, we have 
$
\nabla _w \Phi=0
$
for $x=x'=0$, $w=0$, $\sigma=\pm n(0,0)$. 
According to the curvature property~\eqref{eq1.2}, we have
$
\text{det} ( \nabla_w ^2 \Phi) \neq 0 
$
for  $x=x'=0$, $w=0$, $\sigma=\pm n(0,0)$. 
From the implicit function theorem, there exist $\kappa >0$, such that if 
\begin{equation}\label{under}
|\sigma - n(0, 0)| \leq \kappa\quad {\rm or }\quad |\sigma + n(0, 0)| \leq \kappa
\end{equation}
then the phase $\Phi(x,x', \sigma, w)$ has a unique non-degenerate critical point $w(x,x', \sigma)$ and, by stationary phase, 
under the assumption (\ref{under}), the kernel (\ref{eq.kernel}) is bounded by
$
C\,(1+ \lambda)^{-\frac {d-1 }{2}}.
$
Let us next assume that 
\begin{equation}\label{rakia}
|\sigma - n(0, 0)| >\kappa\quad {\rm and }\quad |\sigma +n(0, 0)| >\kappa\, .
\end{equation}
Then for $w$ close to $0$ and $|x|$ small enough, we obtain by continuity 
\begin{equation}\label{eq.noncrit}
|\sigma - n(x,w)| >\kappa/2 \quad {\rm and }\quad |\sigma + n(x,w)| >\kappa/2\, .
\end{equation}
The kernel of $\nabla_x \nabla_w \varphi_r(x, w)$ is one dimensional and spanned by $n(x, w)$. 
Coming back to the definition of $\Phi$, we deduce that~\eqref{eq.noncrit} implies (for $|x'|$ small enough)
$$
|\nabla_w \Phi (x, x', \sigma,w)| \geq c >0\, .
$$
Consequently, integrating by parts in~\eqref{eq.kernel}, we obtain that under the assumption
(\ref{rakia}) the kernel (\ref{eq.kernel}) is bounded by
$
C_N\, (1+ \lambda)^{-N}
$
which is even better than needed.
This completes the proof of Lemma~\ref{lem.dispersion}.
\end{proof}
The second property of the phase we need is the following:
\begin{lemme}\label{curvature}
Let  $x= (t,z)\in \R\times\R^{d-1}$ where $t=x_1$ and $z=(x_2,\dots,x_{d})$. Then for every
$$
\underline{\omega}=(\underline{\omega_1}, \dots, \underline{\omega_d})\in S^{d-1}
$$
with $\underline{\omega_1}\neq 0$ there exist a neighborhood  $U\subset S^{d-1}$ of $\underline{\omega}$, $\varepsilon>0$
and $\delta>0$ such that, for $\varepsilon/C \leq r \leq C\varepsilon$ and $|x|<\delta$,
the phase $\varphi_{r}(t,z,w)$, where $w\in \R^{d-1}$ is a local coordinate in $U$, is uniformly non degenerate with respect to $(z,w)$.
More precisely
\begin{equation}\label{eq1.0} 
\left|\det_{i,j}\left(\frac{\partial^2 \varphi_{r}(t,z,w)} { \partial z_{j}\partial w_{i}}\right)\right|\geq c >0.
\end{equation}
\end{lemme}
\begin{proof}
Since (\ref{eq1.0}) is independent of the choice of coordinates $w$, it suffices 
to prove it for a particular choice of a coordinate system near $\underline{\omega}$.

For $\omega=(\omega_1,\omega_2,\dots, \omega_d)\in S^{d-1}$ in a small neighborhood of $\underline{\omega}$, we choose $w$ as
$$
w=(w_1,\dots,w_{d-1}):=(\omega_2,\dots,\omega_{d})
$$
which is a coordinate system thanks to the assumption $\underline{\omega_1}\neq 0$. We can also assume that at the point $(t=0, z=0)$, 
the metric is diagonal $g_{i,j}= \delta_{i,j}$. 
Using  Lemma \ref{pak}, we get
\begin{equation}\label{eq1.0-bis} 
\left.\det_{i,j}\left(\frac{\partial^2 \varphi_{r}(t,z,w)} { \partial z_{j}\partial w_{i}}\right)\right|_{(t,z,w)=(0,0,w)}=1.
\end{equation}
We now obtain (\ref{eq1.0})  from (\ref{eq1.0-bis}) by continuity.
\end{proof} 
We next state a corollary of Lemma \ref{curvature}.
\begin{lemme}\label{curvature-bis}
Let $\underline{\omega}^{(1)},\dots, \underline{\omega}^{(N)}$ be $N$ points on $S^{d-1}$. 
Then there exists a splitting of the variable $x= (t,z)\in \R\times\R^{d-1}$ and 
neighborhoods  $U_{j}\subset S^{d-1}$, $j=1,\dots,N$ of $\underline{\omega}^{(j)}$ such
that $\varphi_{r}(t,z,w)$ satisfies (\ref{eq1.0}), where $w$ is a coordinate
in $\cup _{j=1}^{N} U_{j}$.
\end{lemme}
\begin{proof}
Obviously,  there exists a unit vector $e$ such that
$$
e\,\cdot\, \underline{\omega}^{(j)}\neq 0,\qquad j=1,\dots,N.
$$
By performing a rotation, we can assume that $e=(1,0,\dots,0)$ and consequently it suffices to apply  Lemma \ref{curvature}.
\end{proof}
\subsection{Linear estimates}
The dispersion inequality of Lemma~\ref{lem.dispersion} leads to the following estimate.
\begin{lemme}\label{linfty}
Let $(t,z) \in \mathbb{R} \times\mathbb{R}^{d-1}$ be any local system of coordinate near $(0,0)$. Then  the operator
$$
g\in L^2_{w} \longmapsto (T_{\nu}^{r})g(t,z)\in L^{2}( \R_{t}; L^{\infty}( \R^{d-1}_{z}))
$$
is continuous with norm bounded by $C \Lambda(d, \nu) \nu^{-(d-1)/2}$.
\end{lemme}
\begin{proof}
Recall that
$$
(T_{\nu}^{r}f)(t,z)=\int e^{i\nu\varphi_{r}(t,z,w)} a_r(t,z,w,\nu)f(w)dw\, .
$$
Let consider the formal adjoint of $T_{\nu}^{r}$ defined as
$$
(T_{\nu}^{r})^{*}(g)(w):=\int e^{-i\nu\varphi_{r}(t',z',w)} \overline{a_r(t',z',w,\nu)}g(t',z')dt'dz'
$$
According to the classical duality argument which reduces the study of $T_{\nu}^{r}$ to the study
of $T_{\nu}^{r}(T_{\nu}^{r})^{*}$, it is sufficient
to show that the norm of the operator
$$
T_{\nu}^{r}(T_{\nu}^{r})^* \,:\, L^{2}_{t}L^{1}_{z} \longmapsto L^{2}_{t}L^{\infty}_{z}
$$
is bounded by $C[\Lambda(d, \nu) \nu^{-(d-1)/2}]^{2}$. 
But according to Lemma~\ref{lem.dispersion}, the kernel of this operator satisfies~\eqref{eq.dispersion} and as a consequence, 
there exists $C>0$ such that for every $\nu\geq 1$,
\begin{equation}\label{eq.dispfble}
|K(t,z,t',z')|\leq \frac{C}{(1+\nu|t-t'|)^{\frac{d-1}{2}}}\,.
\end{equation}
Using (\ref{eq.dispfble}) and the Young inequality, we get
\begin{equation*}
\|T_{\nu}^{r}(T_{\nu}^{r})^*\, g\|_{L^{2}_{t}L^{\infty}_{z}} 
\leq C \int_{|s|\leq c} \frac {ds} {(1+\nu |s|)^{ \frac{d-1}{2}}}\,\, \|g\|_{L^{2}_{t}L^{1}_{z}} 
\end{equation*}
But clearly
\begin{equation*}
\int_{|s|\leq c} \frac {ds} {(1+\nu |s|)^{ \frac{d-1}{2}}}
\leq {\begin{cases} 
C \nu^{-1/2} & \text{ if } d=2\\
C \nu^{-1}\log (\nu) & \text{ if } d=3\\
C \nu^{-1} & \text{ if } d\geq 4.
\end{cases}}
\end{equation*}
It remains to observe that the right hand-side of the above inequality is equal to $C[\Lambda(d, \nu) \nu^{-(d-1)/2}]^{2}$
which completes the proof of Lemma \ref{linfty}.
\end{proof}
In two space dimensions, we shall need the following extension of Lemma~\ref{linfty}.
\begin{lemme}\label{linfty-bis}
Let $d=2$ and $(t,z) \in \mathbb{R} \times\mathbb{R}^{d-1}$ be any local system of coordinate near $(0,0)$. The operator
$$
g\in L^2_{w} \longmapsto (T_{\nu}^{r})g(t,z)\in L^{4}( \R_{t}; L^{\infty}( \R_{z}))
$$
is continuous with norm bounded by $C\nu^{-1/4}$.
\end{lemme}
\begin{proof}
Similarly to the proof of Lemma~\ref{linfty}, it is sufficient to show that the norm of the operator
$$
T_{\nu}^{r}(T_{\nu}^{r})^* \,:\, L^{4/3}_{t}L^{1}_{z} \longmapsto L^{4}_{t}L^{\infty}_{z}
$$
is bounded by $C\nu^{-1/2}$. The kernel $K(t,z,t',z')$ 
of $T_{\nu}^{r}(T_{\nu}^{r})^*$ satisfies the bound (\ref{eq.dispfble}) with $d=2$.
From the Hardy-Littlewood inequality,
\begin{equation}\label{HL}
\left\|
\int_{-\infty}^{\infty}
\frac{f(t')dt'}{(1+\nu|t-t'|)^{\frac{1}{2}}}
\right\|_{L^{4}(\R_{t})}
\leq C\nu^{-1/2}\|f\|_{L^{4/3}(\R)}\, .
\end{equation}
Therefore
\begin{equation*}
\|T_{\nu}^{r}(T_{\nu}^{r})^*\, g\|_{L^{4}_{t}L^{\infty}_{z}} 
\leq C \nu^{-1/2}\|g\|_{L^{4/3}_{t}L^{1}_{z}} 
\end{equation*}
which completes the proof of Lemma~\ref{linfty-bis}.
\end{proof} 
In the proof of (\ref{tri}), we need the following extension of Lemma~\ref{linfty} for $d\geq 3$.
\begin{lemme}\label{linfty-tris}
Let $d\geq 3$, $p>2$ and $(t,z) \in \mathbb{R} \times\mathbb{R}^{d-1}$ be any local system of coordinate near $(0,0)$. The operator
$$
g\in L^2_{w} \longmapsto (T_{\nu}^{r})g(t,z)\in L^{p}( \R_{t}; L^{\infty}( \R_{z}^{d-1}))
$$
is continuous with norm bounded by $C\nu^{-1/p}$.
\end{lemme}
\begin{proof}
Let $p'$ be such that $\frac{1}{p}+\frac{1}{p'}=1$. It is sufficient to show that the norm of the operator
$$
T_{\nu}^{r}(T_{\nu}^{r})^* \,:\, L^{p'}_{t}L^{1}_{z} \longmapsto L^{p}_{t}L^{\infty}_{z}
$$
is bounded by $C\nu^{-2/p}$. Since for $p>2$,
$$
\left\|\frac{1}{(1+\nu|t|)^{\frac{d-1}{2}}}\right\|_{L^{\frac{p}{2}}(\R_{t})}
\leq C\nu^{-2/p},
$$
using the Young inequality, we get the bound
\begin{equation*}
\left\|\int_{-\infty}^{\infty}\frac{f(t')dt'}{(1+\nu|t-t'|)^{\frac{d-1}{2}}}\right\|_{L^{p}(\R_{t})}\leq C\nu^{-2/p}\|f\|_{L^{p'}(\R)}
\end{equation*}
which completes the proof of Lemma~\ref{linfty-tris} thanks to the bound  (\ref{eq.dispfble}) on the kernel of
$T_{\nu}^{r}(T_{\nu}^{r})^*$.
\end{proof}
\begin{remarque}
Notice that for $p=2$, the proof above still works in dimensions $d\geq 4$.
In the case $p=2$, $d=3$, we have the same difficulty as in the case of the end point Strichartz
estimates on $\R^2$ (see Remark~\ref{kochtat} below).
\end{remarque}
A consequence of Lemma \ref{curvature} is the following statement.
\begin{lemme}\label{l2} Under the assumptions of Lemma~\ref{curvature}, the operator
$$
g\in L^2_{w} \longmapsto (T_{\nu}^{r}g)(t,z)\in L^\infty( \R_{t}; L^2( \R^{d-1}_{z}))
$$
is continuous with norm bounded by $C \nu^{ -(d-1)/2}$.
\end{lemme}
\begin{proof}
In view of  (\ref{eq1.0}), the statement of Lemma \ref{l2}, which can be understood as a refinement of the 
$L^2$ boundedness of the spectral projector, is an immediate consequence of the following generalization of
Plancherel's identity, which we borrow from \cite{Ho1}.
\begin{lemmecite}[non degenerate phase lemma]
Let us consider $\varphi \in C^{\infty}(\R^{n}_{z}\times\R^{n}_{w})$  
and $a\in C_{0}^{\infty}(\R^{n}_{z}\times\R^{n}_{w})$ such that
\begin{eqnarray}\label{ND}
(z,w)\in {\rm supp}(a)\Longrightarrow 
\det\left[\frac{\partial^{2}\varphi}{\partial z\partial w}(z,w)\right]\neq 0.
\end{eqnarray}
There exists $C>0$ such that for every $\lambda\geq 1$, the operator $T_{\lambda}$ 
\begin{equation}\label{def-w}
T_{\lambda}f(z)=\int_{\R^n}e^{i\lambda\varphi(z,w)}a(z,w)f(w)dw
\end{equation}
satisfies,
$$
\|T_{\lambda}f\|_{L^{2}(\R^{n}_{z})}\leq C\lambda^{-\frac{n}{2}}\|f\|_{L^2(\R^n_{w})}\, .
$$
\end{lemmecite}
\end{proof}
\subsection{Multilinear estimates}
Let us first prove (\ref{edno}). We can write
\begin{multline*}
({T}_{\lambda}^{r}f\,{T}_{\mu}^{q}g)(x)\\
 = 
\int_{S^{d-1}}
\int_{S^{d-1}}
e^{i\lambda\varphi_{r}(x,\omega)+i\mu\varphi_{q}(x,\omega')}
a_r(x,\omega,\lambda){a}_q(x,\omega',\mu){f}(\omega){g}(\omega')d\omega d\omega'.
\end{multline*}
We need to evaluate the above expression in $L^2_{x}$. After a partition of unity, we can suppose that 
on the support of 
$$
a_r(x,\omega,\lambda) a_q(x,\omega', \mu),
$$
$(\omega,\omega')$ is close to a fixed point $(\underline{\omega}^{(1)},\underline{\omega}^{(2)})$.
We can therefore use the splitting $x=(t,z)$ of Lemma~\ref{curvature-bis} with $N=2$. 
Using H\"older's inequality, Lemma~\ref{l2} and Lemma \ref{linfty}, we infer
\begin{equation*}
\|{ {T}}_{\lambda}^{r}f\,{{T}}_{\mu}^{q}g\|_{L^2_{t}L^{2}_{z}}
\leq
\|{ {T}}_{\lambda}^{r}f\|_{L^{2}_{t}L^{\infty}_{z}}
\|{{T}}_{\mu}^{q}g\|_{L^{\infty}_{t}L^{2}_{z}}
\leq
C\Lambda(d,\lambda)(\lambda \mu)^{-\frac{d-1}{2}}\|f\|_{L^2_{w}}\|g\|_{L^2_{w}}.
\end{equation*}
This completes the proof of (\ref{edno}).

We next prove (\ref{dve}). Let us write
\begin{multline*}
({T}_{\lambda}^{r}f\,{T}_{\mu}^{q}g\, {T}_{\nu}^{s}h )(x)
 = 
\int_{S^{1}}\int_{S^{1}}\int_{S^{1}}
e^{i\lambda\varphi_{r}(x,\omega)+i\mu\varphi_{q}(x,\omega')+i\nu\varphi_{s}(x,\omega'')}
\\
a_r(x,\omega,\lambda){a}_q(x,\omega',\mu){a}_s(x,\omega'',\nu)
{f}(\omega){g}(\omega'){h}(\omega'')
d\omega d\omega'd\omega''\, .
\end{multline*}
After a partition of unity we can assume that  $(\omega'',\omega',\omega)$
is close to a fixed point $(\underline{\omega}^{(1)},\underline{\omega}^{(2)},\underline{\omega}^{(3)})$.
With the splitting $x=(t,z)$ of Lemma~\ref{curvature-bis} with $N=3$, using 
Lemma \ref{l2}, Lemma \ref{linfty-bis} and the H\"older inequality, we get 
\begin{multline*}
\|{ {T}}_{\lambda}^{r}f\,{{T}}_{\mu}^{q}g\,{T}_{\nu}^{s}h\|_{L^{2}_{t}L^{2}_{z}}
\leq 
\|{ {T}}_{\lambda}^{r}f\|_{L^{4}_{t}L^{\infty}_{z}}
\|{{T}}_{\mu}^{q}g\|_{L^{4}_{t}L^{\infty}_{z}}
\|{T}_{\nu}^{s}h\|_{L^{\infty}_{t}L^{2}_{z}}
\\
\leq
C
\lambda^{-\frac 1 4}
\mu^{-\frac 1 4}
{\nu^{-\frac{1} 2}}
\|f\|_{L^2}\|g\|_{L^2}\|h\|_{L^2}\, .
\end{multline*}
This completes the proof of (\ref{dve}).

We finally prove (\ref{tri}). We can again use the splitting $x=(t,z)$ of 
Lemma~\ref{curvature-bis} with $N=3$. For $p\gg 2$ and $q>2$ such that $\frac{1}{p}+\frac{1}{q}=\frac{1}{2}$, 
a use of Lemma~\ref{linfty-tris} gives the bound
\begin{multline*}
\|{ {T}}_{\lambda}^{r}f\,{{T}}_{\mu}^{q}g\,{T}_{\nu}^{s}h\|_{L^{2}_{t}L^{2}_{z}}
\leq 
\|{ {T}}_{\lambda}^{r}f\|_{L^{p}_{t}L^{\infty}_{z}}
\|{{T}}_{\mu}^{q}g\|_{L^{q}_{t}L^{\infty}_{z}}
\|{T}_{\nu}^{s}h\|_{L^{\infty}_{t}L^{2}_{z}}
\\
\leq
C
\lambda^{-\eta}\mu^{-\frac{1}{2}+\eta}\nu^{-1}
\|f\|_{L^2}\|g\|_{L^2}\|h\|_{L^2}\, ,
\end{multline*}
where $\eta=\frac{1}{p}$. This ends the proof (\ref{tri}) and completes the proof of  Theorem \ref{thm2}. 
\qed
\begin{remarque}\label{kochtat}
As pointed to us by Koch and Tataru~\cite{KT04}, another approach to these multilinear estimates would be, after a suitable 
micro-localization, to particularize one variable ($t$ in the exposition above) and see the equation satisfied by the approximated 
spectral projector
$$ 
({\mathbf \Delta }+\lambda^ 2 )\chi_\lambda (u) = \mathcal{O}_{L^2} (\lambda)
$$
as a {\em semi-classical} evolution equation of the type 
$$ 
(ih \partial_t +Q(t,z,hD_z)) \chi_\lambda (u)= \mathcal{O}_{L^2} (h), \qquad h = \lambda^{-1}\, .
$$
Then Lemmas~\ref{linfty}, \ref{linfty-bis} are simply the (semi-classical) Strichartz estimates which can be proved by 
using the approach in~\cite{BGT1}.
\end{remarque}
\section{Preliminaries to the proof of Theorem \ref{thm1}}\label{sec.3}
In this section $(M,g)$ is an arbitrary Riemannian manifold of dimension $d$.
Our first purpose is to introduce the basic localization operators $\Delta _{N}$ and $\Delta _{NL}$ which are naturally
related to the Sobolev spaces and the Bourgain spaces on $M$ respectively. We establish some basic bounds related to
$\Delta _{N}$ and $\Delta _{NL}$. The main purpose of this section is
to show that the well-posedness of the Cauchy problem (\ref{intr-1}) in 
the Sobolev space $H^s(M)$ is a consequence of nonlinear estimates in the Bourgain 
spaces associated to the Laplace operator ${\mathbf \Delta}$. 
This reduction is now classical (see e.g. \cite{Gi}).
\subsection{Bourgain spaces and basic localization operators}
Since $M$ is compact, ${\mathbf \Delta}$ has a compact resolvent and thus the spectrum of ${\mathbf \Delta}$ is discrete. 
Let $e_{k}\in L^2(M)$, $k\in\N$ be an orthonormal basis of eigenfunctions of ${\mathbf \Delta}$ associated to eigenvalues $\lambda_k$.
Denote by $P_k$ the orthogonal projector on $e_{k}$. The Sobolev space $H^{s}(M)$ is therefore equipped with the norm (with $\langle x \rangle = (1 + |x|^2)^{1/2}$),
$$
\|u\|_{H^s(M)}^2= \sum_k \langle \lambda_{k}\rangle^s \|P_k u\|_{L^2(M)}^2.
$$ 
The Bourgain space (or conormal Sobolev space) $X^{s,b}(\R\times M)$ is equipped with the  norm
\begin{equation*}
\|u\|_{{X}^{s,b}(\R\times M)}^2=\sum_{k}\langle \lambda_{k}\rangle^{s }\|\langle\tau+\lambda_{k}\rangle^{b}
\widehat{P_k u}(\tau)\|^2_{L^{2}(\R_{\tau}\times M)}
= \|e^{it {\mathbf \Delta}} u(t, \cdot)\|^2_{H^b( {\R}_{t}\,;\,H^s(M))},
\end{equation*}
where $\widehat{P_k u}(\tau)$ denotes the Fourier transform of $P_k u$ with respect to the time variable.

Let us first recall that for $b>1/2$  the space ${X}^{s,b}(\R\times M)$ is continuously embedded in $C(\R\,;\,H^{s}(M))$.
For $u\in C_{0}^{\infty}(\R\times M)$, we write 
$$
P_k u(t)=\frac{1}{2\pi}\int_{-\infty}^{\infty}\frac{\langle\tau+\lambda_k\rangle^{b}}{\langle\tau+\lambda_k\rangle^{b}}
\widehat{P_k u}(\tau)e^{it\tau} d\tau.
$$
For $b>\frac{1}{2}$, we get by the Cauchy-Schwarz inequality, applied in $\tau$,
\begin{equation}\label{eq2.400}
\langle\lambda_{k}\rangle^{s/2}|P_k u(t)|\leq C\left\{\int_{-\infty}^{\infty}
\langle\lambda_{k}\rangle^{s}
\langle\tau+\lambda_k\rangle^{2b}|\widehat{P_k u}(\tau)|^{2}d\tau\right\}^{\frac{1}{2}}.
\end{equation}
Squaring (\ref{eq2.400}), integrating over $M$ and summing over $k$ yields,
\begin{equation}\label{Sobolev}
\|u\|_{L^{\infty}(\R\,;\,H^{s}(M))}\leq C\|u\|_{X^{s,b}(\R\times M)},\quad b>\frac{1}{2}.
\end{equation}
For $u\in C^{\infty}(M)$ and $N\geq 1$, we define the projector $\Delta _{N}$ as
$$
\Delta _{N}(u):=\sum_{k\,:\, N\leq \langle\lambda_k\rangle^{\frac{1}{2}}< 2N}P_{k}u\, .
$$
We now state a basic bound for $\Delta _{N}$.
\begin{lemme}\label{sob-g}
There exists a constant $C$ such that for every $q\in[2,\infty]$, every $u\in L^{2}(M)$,
$$
\|\Delta _{N}(u)\|_{L^{q}(M)}\leq C\,N^{\frac{d}{2}-\frac{d}{q}}\|\Delta _{N}(u)\|_{L^{2}(M)}\, .
$$
\end{lemme}
\begin{proof}
The assertion clearly holds for $q=2$. We next prove it for $q=\infty$.
Let us write
$
\Delta _{N}=\sum_{j=0}^{N-1}\Delta _{N,j},
$
where
\begin{equation}\label{cut}
\Delta _{N,j}:=\sum_{k\,:\, N+j\leq \langle\lambda_k\rangle^{\frac{1}{2}}< N+j+1}P_{k}\, .
\end{equation}
Due to the  Weyl estimate (see Lemma \ref{sogge})
$$
\|\Delta _{N,j}(u)\|_{L^{\infty}(M)}\leq CN^{\frac{d-1}{2}}\|\Delta _{N,j}(u)\|_{L^{2}(M)}
$$
and due to the triangle and Cauchy-Schwarz inequalities
\begin{multline*}
\|\Delta _{N}(u)\|_{L^{\infty}(M)}\\
\leq CN^{\frac{d-1}{2}}\Big(\sum_{j=0}^{N-1}\|\Delta _{N,j}(u)\|^{2}_{L^{2}(M)}\Big)^{\frac{1}{2}}
\Big(\sum_{j=0}^{N-1}1^{2}\Big)^{\frac{1}{2}}
= CN^{\frac{d}{2}}\|\Delta _{N}(u)\|_{L^{2}(M)}\, .
\end{multline*}
By H\"older's inequality, this completes the proof of Lemma \ref{sob-g}.
\end{proof}
For $N\geq 1$ a dyadic integer, i.e. $N=2^{n}$, $n\in \N$, we define the operator $S_N$ as
$$
S_{N}:=\sum_{N_1\leq N}\Delta _{N_{1}},
$$ 
where the sum is taken over all dyadic integers $N_1$ smaller or equal to $N$. We also define $S_{\frac{1}{2}}$ by
$S_{\frac{1}{2}}(u):=0$.

Next, for $u\in C_{0}^{\infty}(\R\times M)$ and $N,L$ positive integers, we define the localization operators $\Delta _{NL}$ as
$$
\Delta _{NL}(u):=\frac{1}{2\pi}\sum_{k\,:\, N\leq \langle\lambda_k\rangle^{\frac{1}{2}}< 2N}
\int_{L\leq\langle\tau+\lambda_k\rangle\leq 2L}\widehat{P_k u}(\tau)e^{it\tau}d\tau\, .
$$
It is easy to check that  $\Delta _{NL}$ is a projector.
It follows from this definition that for every $s$, $b$ there exists $C>0$ such that
\begin{equation}\label{eq-bis}
\frac{1}{C}\|\Delta _{NL}(u)\|_{X^{s,b}(\R\times M)}\leq L^{b}N^{s} \|\Delta _{NL}(u)\|_{L^2(\R\times M)}\leq 
C\|\Delta _{NL}(u)\|_{X^{s,b}(\R\times M)}\, ,
\end{equation}
and 
\begin{equation}\label{harcterizatzia}
\frac{1}{C}\sum_{N,L}L^{2b}N^{2s} \|\Delta _{NL}(u)\|_{L^2}^{2}\leq\|u\|_{X^{s,b}}^{2}
\leq C\sum_{N,L}L^{2b}N^{2s} \|\Delta _{NL}(u)\|_{L^2}^{2}\, ,
\end{equation}
where the sums is taken over all dyadic values of $N$ and $L$, i.e. $N=2^{n}$, $L=2^{l}$, $n,l\in \N$.
We now state a basic bound for $\Delta _{NL}$.
\begin{lemme}\label{sob-co-g}
There exists a constant $C$ such that for every $p,q\in[2,\infty]$, every $u\in L^{2}(\R\times M)$,
$$
\|\Delta _{NL}(u)\|_{L^{p}(\R\,;\,L^{q}(M))}
\leq C
\,L^{\frac{1}{2}-\frac{1}{p}}
N^{\frac{d}{2}-\frac{d}{q}}\|\Delta _{NL}(u)\|_{L^{2}(\R\times M)}\, .
$$
\end{lemme}
\begin{proof}
Since $\Delta _{N}\Delta _{NL}=\Delta _{NL}$, a use of Lemma \ref{sob-g} yields
$$
\|\Delta _{NL}(u)\|_{L^{p}(\R\,;\,L^{q}(M))}
\leq CN^{\frac{d}{2}-\frac{d}{q}}\|\Delta _{NL}(u)\|_{L^{p}(\R\,;\,L^{2}(M))}.
$$
Therefore, we only need to consider the case $q=2$. 
Using that $\widehat{P_k u}(\tau)=P_{k}\widehat{u}(\tau)$, we can write
$$
\|\Delta _{NL}(u)(t)\|_{L^{2}(M)}^{2}\leq C
\sum_{k\,:\, N\leq \langle\lambda_k\rangle^{\frac{1}{2}}< 2N}\,\int_{M}
\left|
\int_{L\leq\langle\tau+\lambda_k\rangle\leq 2L}\widehat{P_k u}(\tau)e^{it\tau}d\tau
\right|^{2}\, .
$$
Since the integration over $\tau$ is on a region of size $L$, using the Cauchy-Schwarz inequality in
$\tau$ and the Plancherel identity yields
\begin{equation*}
\|\Delta _{NL}(u)(t)\|_{L^{2}(M)}^{2}\leq CL \sum_{k\,:\, N\leq \langle\lambda_k\rangle^{\frac{1}{2}}< 2N}\,
\,\int_{M}\int_{-\infty}^{\infty}|\widehat{P_k u}(\tau)|^{2}d\tau\leq CL\|u\|_{L^{2}(\R\times M)}^{2}\,\, .
\end{equation*}
Applying the last inequality to $\Delta _{NL}(u)$ instead of $u$
and using that $\Delta _{NL}$ is a projector gives
$$
\|\Delta _{NL}(u)\|_{L^{\infty}(\R\,;\,L^{2}(M))}
\leq C
\,L^{\frac{1}{2}}\|\Delta _{NL}(u)\|_{L^{2}(\R\times M)}\, .
$$
The assertion of the lemma trivially holds for $p=q=2$ and therefore the proof of Lemma \ref{sob-co-g} is completed by
H\"older's inequality.
\end{proof}
\subsection{Reduction to nonlinear estimates in Bourgain spaces}
The starting point is to consider the integral equation (Duhamel form)
\begin{equation}\label{Duhamel}
u(t)=e^{it{\mathbf \Delta}}u_{0}- i\int_{0}^{t}e^{i(t-\tau){\mathbf \Delta}}F(u(\tau))d\tau\, .
\end{equation}
At least for classical (smooth) solutions $u$ the integral equation (\ref{Duhamel}) is equivalent to  (\ref{intr-1}).
For that reason we solve (\ref{Duhamel}) by the Picard iteration in a suitable functional setting and thus we get solutions
of (\ref{Duhamel}). Notice that this is achieved classically if $s>3/2$ by taking $X= C([-T, T]; H^s(M))$. Therefore we shall restrict the study to the case $s\in [1, 3/2]$. In this case of low regularity solutions of (\ref{Duhamel}) the information we have for $u$ and $F(u)$ should be
strong enough to conclude that we get solutions of (\ref{intr-1}) too.

For $T>0$, we define the restriction space $X^{s,b}_{T}:=X^{s,b}([-T,T]\times M)$, equipped with the norm
$$
\|u\|_{X^{s,b}_{T}}=
\inf\{\|w\|_{X^{s,b}},\quad w\in X^{s,b} \quad{\rm with }\quad w|_{[-T,T]}=u\}.
$$
For $b>1/2$, the space $X^{s,b}_{T}$ is continuously embedded in $C([-T,T]\,;\,H^{s}(M))$ (see (\ref{Sobolev})) 
and $X^{s,b}_{T}$ will be the
space for the solutions of (\ref{Duhamel}) on $[-T,T]$. 
The next proposition contains the basic reduction to nonlinear estimates in $X^{s,b}$.
\begin{proposition}\label{basic-reduction}
Assume that there exists $(b,b')\in\R^2$ satisfying
\begin{equation}\label{restrictions}
0<b'<\frac{1}{2}<b,\quad b+b'<1
\end{equation}
such that for every $s\in [1,2)$ there exists a constant $C$ and $p$ such that for every $u\in X^{s,b}$,
\begin{equation}\label{parvo-red}
\|F(u)\|_{X^{s,-b'}(\R\times M)}\leq C\Big(1+\|u\|_{X^{1,b}(\R\times M)}^{p}\Big)\|u\|_{X^{s,b}(\R\times M)}\, ,
\end{equation}
and for every $u,v\in X^{s,b}$,
\begin{multline}\label{vtoro-red}
\|F(u)-F(v)\|_{X^{s,-b'}(\R\times M)}
\leq \\
C\Big(1+\|u\|_{X^{s,b}(\R\times M)}^{p}+\|v\|_{X^{s,b}(\R\times M)}^{p}\Big)\|u-v\|_{X^{s,b}(\R\times M)}\, .
\end{multline}
Then
\begin{enumerate}
\item\label{item-1}
For every bounded subset $B$ of $H^1(M)$ there exists $T>0$ such that if $u_0\in B$ then there exists a unique solution $u$ of
(\ref{Duhamel}) on $[-T,T]$ in the class $X^{1,b}_{T}$. 
Moreover the map  $u_0\mapsto u\in C([-T,T]\,;\, H^{1}(M))$ is  Lipschitz continuous on $B$.
\item\label{item-2}
If in addition $u_{0}\in H^{s}(M)$ then $u\in C([-T,T]\,;\, H^{s}(M))$.
\item
The function $u$ is a solution of (\ref{intr-1}) in the distributional sense.
\item
If in addition we suppose that $V(z)\geq -C(1+|z|)^{\beta}$, $\beta<2+4/d$
then the results above hold for any arbitrarily large $T$.
\item
For $u_{0}\in H^{s}(M)$, $s>3/2$ the solution is unique in $C([-T,T]\,;\, H^{s}(M))$.
\end{enumerate}
\end{proposition}
\begin{proof}
Let $\psi\in C_{0}^{\infty}(\R)$ be equal to $1$ on $[-1,1]$. The identity
$$
\|\psi(t)e^{it{\mathbf \Delta}}\,u_{0}\|_{X^{s,b}(\R\times M)}=\|\psi\|_{H^{b}(\R)}\|u_0\|_{H^s(M)}
$$
follows from the definition of $X^{s,b}(\R\times M)$ and therefore for $T\leq 1$
\begin{equation}\label{free}
\|e^{it{\mathbf \Delta}}\,u_{0}\|_{X^{s,b}_T}\leq C\|u_0\|_{H^s(M)}\, .
\end{equation}
The bound
\begin{equation}\label{eq2.17}
\Big
\|\psi(t/T)\int_{0}^{t}e^{i(t-\tau){\mathbf \Delta}}F(\tau)d\tau
\Big\|_{X^{s,b}(\R_t\times M)}
\leq C T^{1-b-b'}\|F\|_{X^{s,-b'}(\R\times M)},
\end{equation}
holds for $T\leq 1$ and $(b,b')$ satisfying (\ref{restrictions}). 
Indeed (see \cite[Proposition 2.11]{BGT'03}) estimate (\ref{eq2.17}) follows from the {\it one dimensional} inequality
\begin{eqnarray}\label{basic}
\|\psi(t/T)\int_{0}^{t}g(\tau)d\tau\|_{H^{b}(\R)}\leq C T^{1-b-b'}\|g\|_{H^{-b'}(\R)}.
\end{eqnarray}
A proof of (\ref{basic}) can be found in \cite{Gi}.

Using (\ref{eq2.17}) and the assumptions of the proposition we obtain the estimates
\begin{equation}\label{tame}
\Big\|
\int_{0}^{t}e^{i(t-\tau){\mathbf \Delta}}F(u(\tau))d\tau
\Big\|_{X^{s,b}_{T}}
\leq C T^{1-b-b'}\big(1+\|u\|_{X^{1,b}_{T}}^{p}\big)\|u\|_{X^{s,b}_{T}}
\end{equation}
and
\begin{multline}\label{diff}
\Big\|
\int_{0}^{t}e^{i(t-\tau){\mathbf \Delta}}(F(u(\tau))-F(v(\tau)))d\tau
\Big\|_{X^{s,b}_{T}}
\leq \\
C T^{1-b-b'}
\big(1+\|u\|_{X^{s,b}_{T}}^{p}+\|v\|_{X^{s,b}_{T}}^{p}\big)\|u-v\|_{X^{s,b}_{T}}\, ,
\end{multline}
provided $T\leq 1$ and $(b,b')$ satisfying (\ref{restrictions}). 
Let $B$ be a bounded subset of $H^1(M)$.
It results from  (\ref{free}), (\ref{tame}) and (\ref{diff}) with $s=1$ 
that there exists $T\ll 1$ such that for every $u_{0}\in B$ the right hand-side of (\ref{Duhamel}) is a contraction
in a suitable ball of $X^{1,b}_{T}$ with a unique fixed point which is the solution of (\ref{Duhamel}).
The uniqueness in the class $X^{1,b}_{T}$ and the Lipschitz continuity of the flow map follow from (\ref{diff}).
Suppose now that $u_{0}\in H^{s}(M)$. Then as before it follows from (\ref{free}), (\ref{tame}) and (\ref{diff}) that
we can find $\widetilde{T}\leq T$ such that we can identify $u|_{[-\widetilde{T},\widetilde{T}]}$ as the
unique solution of (\ref{Duhamel}) on $[-\widetilde{T},\widetilde{T}]$
in the class $X^{s,b}_{\widetilde{T}}\subset X^{1,b}_{\widetilde{T}} $. 
In particular $u(t,\cdot)\in H^{s}(M)$ for $t\in [-\widetilde{T},\widetilde{T}]$.
Then by a bootstrap and the tame estimate (\ref{tame})
we conclude that $u(t,\cdot)\in H^{s}(M)$ for $t\in [-T,T]$.
Thanks to (\ref{parvo-red}), we obtain that $F(u(t))\in X^{1,-b'}_{T}$ and since $b'<1/2$, we infer that 
$$
\partial_{t}\Big[\int_{0}^{t}e^{-i\tau{\mathbf \Delta}}F(u(\tau))d\tau\Big]= e^{-it{\mathbf \Delta}}F(u(t))
$$
in the distributional sense
which implies that $u$ is a solution of the original PDE (\ref{intr-1}) in the distributional sense.
If $u_{0}\in H^{2}(M)$ then, thanks to the propagation of the $H^s$ regularity assertion, 
one can take the scalar product of (\ref{intr-1}) with $u$ and $u_t$ and it results that the conservations laws 
(\ref{intr-2}) hold. If $u_{0}\in H^{1}(M)$, we can approximate in $H^1(M)$ the function $u_0$ with a sequence $(u_{0,n})$
such that $u_{0,n}\in H^{2}(M)$. If we denote by $u_{n}(t)$ the corresponding solutions of (\ref{intr-1}), thanks
to the propagation of the regularity we obtain then $u_{n}(t)$ enjoy the conservation laws (\ref{intr-2}) for $t$ on the 
time of existence of $u(t)$. Finally we can pass to the limit $n\rightarrow\infty$ and thanks to the $H^1$ continuity
of the conservation laws functionals, we deduce that $u(t)$ satisfies the conservation laws (\ref{intr-2}).
If we suppose that $V(z)$ satisfies $V(z)\geq -C(1+|z|)^{\beta}$, $\beta<2+4/d$, using the Gagliardo-Nirenberg inequalities, 
we obtain that there exists $\theta\in ]0,2[$
such that
$$
\int_{M}V(u(t))
\geq
-C\|u(t)\|_{L^{2}}^{\beta-\theta}\|u(t)\|_{H^{1}}^{\theta}-C \, .
$$
Therefore, 
the conservation laws (\ref{intr-2})
yield a bound independent with respect to $t$ for $\|u(t)\|_{H^{1}(M)}$ which allows to reiterate the local existence argument and
thus to achieve the existence of $u(t)$ on an arbitrary time interval.
Finally, thanks to the Sobolev embedding $H^{s}(M)\subset L^{\infty}(M)$, $s>3/2$ and the propagation of regularity, we easily obtain that 
if $u_{0}\in H^{s}(M)$, $s>3/2$ then the uniqueness holds in the class $C([-T,T]\,;\, H^{s}(M))$.
This completes the proof of Proposition \ref{basic-reduction}. 
\end{proof}
With Proposition \ref{basic-reduction} in hand the assertion of Theorem \ref{thm1} follows from the following statement.
\begin{theoreme}\label{thm4}
Let $M=S^3$ or $M=S^2_{\rho}\times S^1$ endowed with the standard metrics. 
For every $1<\alpha<5$ and $s\in[1,2)$ there exists $(b,b')\in\R^2$ satisfying (\ref{restrictions}) such that 
(\ref{parvo-red}) and (\ref{vtoro-red}) hold.
\end{theoreme}
The next two sections are devoted to the proof of Theorem \ref{thm4}.
\section{Bilinear Strichartz estimates and applications}
\label{sec.4}
In this section we prove Theorem~\ref{thm4} for $M=S^3$ with the standard metric. 
\subsection{Bilinear Strichartz estimates on $S^3$}
In the case $M=S^3$ the eigenvalues of $-{\mathbf \Delta}$ are $\lambda_{k}=k^{2}-1$, $k\geq 1$ and this fact plays a 
key role in the analysis. The starting point is the following bilinear improvement of the $L^4$
Strichartz inequality on $S^3$ established in \cite{BGT1}.
\begin{proposition}\label{str-bl-s3}
For every interval $I\subset\R$, every $\varepsilon>0$ 
there exists a constant $C$ such that for every $N_1,N_2\geq 1$, every $f_1,f_2\in L^{2}(M)$,
\begin{equation*}
\big\|\prod_{j=1}^{2}e^{it{\mathbf \Delta}}(\Delta _{N_j}f_{j})\big\|_{L^{2}(I\times M)}
\leq C (\min(N_1,N_2))^{\frac{1}{2}+\varepsilon}\prod_{j=1}^{2}\|\Delta _{N_j}f_j\|_{L^{2}(M)}\, .
\end{equation*}
\end{proposition}
\begin{proof}
By a time translation we can suppose that $I=[0,T]$. Moreover,
for $f\in L^{2}(M)$ the function $e^{it{\mathbf \Delta}}f$ is periodic with respect to $t$ and therefore it suffices to give the proof
with $T=2\pi$. 
Let us write
$$
\prod_{j=1}^{2}e^{it{\mathbf \Delta}}(\Delta _{N_j}f_{j})=\sum_{N_{j}\leq \langle\lambda_{k_{j}}\rangle^{\frac{1}{2}}< 2N_j}
e^{-it(\lambda_{k_1}+\lambda_{k_2})}P_{k_1}(f_1)P_{k_2}(f_2)\, .
$$
Using the Parseval identity with respect to $t$ we get
$$
\big\|\prod_{j=1}^{2}e^{it{\mathbf \Delta}}(\Delta _{N_j}f_{j})\big\|_{L^{2}([0,2\pi]\times M)}^{2}
=
\sum_{\tau\in\Z}
\Big\|
\sum_{\tau=\lambda_{k_1}+\lambda_{k_2}} 
P_{k_1}(f_1)P_{k_2}(f_2)
\Big\|_{L^{2}(M)}^{2},
$$
where the summation over $(k_1,k_2)$ is restricted to $N_{j}\leq  \langle\lambda_{k_{j}}\rangle^{\frac{1}{2}}< 2N_j$, $j=1,2$.
Applying the triangle inequality for the $L^{2}(M)$ norm, 
the Cauchy-Schwarz inequality in the summation over $(k_1,k_2)$, 
and the bilinear estimate of Theorem~\ref{thm2} for $d=3$ yields that for every $\varepsilon>0$,
\begin{multline*}
\big\|\prod_{j=1}^{2}e^{it{\mathbf \Delta}}(\Delta _{N_j}f_{j})\big\|_{L^{2}([0,2\pi]\times M)}^{2}
\\\leq 
C_{\varepsilon}\,(\min(N_1,N_2))^{1+\varepsilon}
\sup_{\tau\in\Z}\alpha_{N_1,N_2}(\tau)\prod_{j=1}^{2}\|\Delta _{N_j}f_j\|_{L^{2}(M)}^{2}\, ,
\end{multline*}
where
$$
\alpha_{N_1,N_2}(\tau)=\#
\left\{
(k_1,k_2)\in \N^{2}\,:\, \tau+2=k_{1}^{2}+k_{2}^{2},\, N_{j}\leq  \langle\lambda_{k_{j}}\rangle^{\frac{1}{2}}< 2N_j,\, j=1,2
\right\}.
$$
We claim that $\alpha_{N_1,N_2}(\tau)\leq C_{\varepsilon}N^{\varepsilon}$. Indeed this follows from the next lemma.
\begin{lemme}\label{le3.2}
For every $\varepsilon>0$ there exists $C>0$ such that for every positive integers $\tau$ and $N$,
\begin{equation}\label{chudo}
\#\{(k_1,k_2)\in\mathbb{N}^2\,:\, N\leq k_{1} < 2N,\quad k_1^2+k_2^2=\tau\}\leq C N^{\varepsilon}\, .
\end{equation}
\end{lemme}
\begin{proof}
This lemma already appeared in \cite{BGT'03} (see \cite[Lemma 3.2]{BGT'03}).
We recall the proof. For $\tau\leq 10N^{4}$ it follows from the divisor bound in the ring of 
Gaussian integers which is a Euclidean division domain.
For $\tau\geq 10N^{4}$ there is at most one value of $(k_1,k_2)$ satisfying the imposed restriction since in this case 
$k_2$ should range in an interval of size smaller than one.
Hence for $\tau\geq 10N^{4}$ the left hand-side of (\ref{chudo}) is bounded by $1$.
This completes the proof of Lemma~\ref{le3.2}.
\end{proof}
Proposition \ref{str-bl-s3} now readily follows from Lemma \ref{le3.2}.
\end{proof}
\subsection{Using bilinear Strichartz estimates}
From now on we simply assume that $M$ is a three dimensional compact manifold satisfying Proposition~\ref{str-bl-s3}.
Proceeding as in~\cite[Section 3.2]{BGT'03} one can show, for instance, that three dimensional Zoll manifolds have this property. 
As a consequence it can be remarked that in fact Theorem~\ref{thm1} holds for any such manifold.
\vskip .5cm
First we deduce from Proposition~\ref{str-bl-s3} the following bilinear estimate in the $X^{s,b}$ context.
\begin{proposition}\label{bolen}
For every $\varepsilon>0$ there exist $\beta<\frac{1}{2}$ and $C>0$ such that for every $N_1,N_2,L_1,L_2\geq 1$, 
every $u_1,u_2\in L^{2}(\R\times M)$,
\begin{equation*}
\big\|\prod_{j=1}^{2}\Delta _{N_{j}L_{j}}(u_{j})\big\|_{L^{2}(\R\times M)}
\leq C(L_1 L_2)^{\beta}
(\min(N_1,N_2))^{\frac{1}{2}+\varepsilon}\prod_{j=1}^{2}\|\Delta _{N_j L_j}(u_j)\|_{L^{2}(\R\times M)}\, .
\end{equation*}
\end{proposition}
\begin{proof}
Let us suppose that $N_1\leq N_2$. 
Using Lemma \ref{sob-co-g} and the H\"older inequality we can write
\begin{multline}\label{parva-pista}
\big\|\prod_{j=1}^{2}\Delta _{N_{j}L_{j}}(u_{j})\big\|_{L^{2}(\R\times M)}
\leq
\\
\leq
\|\Delta _{N_{1}L_{1}}(u_{1})\|_{L^{4}(\R\,;\,L^{\infty}(M))}
\|\Delta _{N_{2}L_{2}}(u_{2})\|_{L^{4}(\R\,;\,L^{2}(M))}
\leq
\\
\leq 
C\, N_{1}^{\frac{3}{2}}(L_1L_2)^{\frac{1}{4}}
\prod_{j=1}^{2}\|\Delta _{N_j L_j}(u_j)\|_{L^{2}(\R\times M)}\, .
\end{multline}
Estimate (\ref{parva-pista}) is better than the needed one with respect to the $L_j$ localization but is far from
the needed one with respect to the $N_j$ localization.

We now estimate $\big\|\prod_{j=1}^{2}\Delta _{N_{j}L_{j}}(u_{j})\big\|_{L^{2}}$ by means of Proposition \ref{str-bl-s3}.
It is indeed possible thanks to the following lemma.
\begin{lemme}\label{claim}
For every $b\in ]\frac{1}{2},1]$, every $\delta>0$, there exists $C_{b,\delta}$ such that for every $u_1,u_2\in
X^{0,b}(\R\times M)$, every $1\leq N_1 \leq N_2$,
\begin{equation*}
\big\|\prod_{j=1}^{2}\Delta _{N_{j}}(u_{j})\big\|_{L^{2}(\R\times M)}
\leq C_{b,\delta}\, 
N_1^{\frac{1}{2}+\delta}\prod_{j=1}^{2}\|\Delta _{N_j}(u_j)\|_{X^{0,b}(\R\times M)} \, .
\end{equation*}
\end{lemme}
\begin{proof}
Let us set $v_{j}(t):=e^{it{\mathbf \Delta}}\Delta _{N_{j}}(u_{j})(t)$, $j=1,2.$
Then we can write
$$
\Delta _{N_{j}}(u_{j})(t)=\frac{1}{2\pi}\int_{-\infty}^{\infty}e^{it\tau}\, e^{-it{\mathbf \Delta}}\widehat{v_{j}}(\tau)\,d\tau\, .
$$
Therefore
$$
\prod_{j=1}^{2}\Delta _{N_{j}}(u_{j})(t)=\frac{1}{4\pi^{2}}\int_{-\infty}^{\infty}\int_{-\infty}^{\infty}
e^{it(\tau_1+\tau_2)}
\big(\prod_{j=1}^{2}
e^{-it{\mathbf \Delta}}\widehat{v_{j}}(\tau_j)\big)\,d\tau_1 d\tau_2\, .
$$
Using the triangle inequality and Proposition \ref{str-bl-s3} gives that for every unit interval $I\subset\R$, 
every $\delta>0$ there exists $C_{\delta}$ such that
$$
\big\|\prod_{j=1}^{2}\Delta _{N_{j}}(u_{j})\big\|_{L^{2}(I\times M)}
\leq 
C_{\delta}N_{1}^{\frac{1}{2}+\delta} 
\int_{-\infty}^{\infty}\int_{-\infty}^{\infty}
\prod_{j=1}^{2}\|\widehat{v_{j}}(\tau_j)\|_{L^{2}(M)}\,d\tau_1 d\tau_2\, .
$$
Hence using the Cauchy-Schwarz inequality in $(\tau_1,\tau_2)$ gives for $b>1/2$,
\begin{multline}\label{split1}
\big\|\prod_{j=1}^{2}\Delta _{N_{j}}(u_{j})\big\|_{L^{2}(I\times M)}
\leq 
C_{b,\delta}N_{1}^{\frac{1}{2}+\delta}
\prod_{j=1}^{2}\|\langle\tau\rangle^{b}\widehat{v_{j}}(\tau)\|_{L^{2}(\R_{\tau}\times M)}
=
\\
=
C_{b,\delta}\, 
N_1^{\frac{1}{2}+\delta}\prod_{j=1}^{2}\|\Delta _{N_j}(u_j)\|_{X^{0,b}(\R\times M)} \, .
\end{multline}
Using a partition of unity, we can find a $\psi\in C_{0}^{\infty}(\R)$, supported in $[0,1]$
such that
\begin{equation}\label{split}
\Delta _{N_{j}}(u_{j})(t)=\sum_{n\in\Z}\psi\big(t-\frac{n}{2}\big)\Delta _{N_{j}}(u_{j}(t))
=
\sum_{n\in\Z}\Delta _{N_{j}}\Big(\psi\big(t-\frac{n}{2}\big)u_{j}(t)\Big)\, .
\end{equation}
Notice that if for $u\in X^{0,b}(\R\times M)$, $b\in [0,1]$, we set
$u_{n}(t)=\psi\big(t-\frac{n}{2}\big)u(t)$ then
\begin{equation}\label{split3}
\sum_{n\in\Z}\|u_{n}\|_{X^{0,b}(\R\times M)}^{2}\leq C\|u\|_{X^{0,b}(\R\times M)}^{2}\, .
\end{equation}
Indeed (\ref{split3}) is straightforward for $b=0$ and $b=1$ and it follows by complex interpolation for $b\in ]0,1[$.
Using the almost disjointness of the supports of $\psi\big(t-\frac{n}{2}\big)$, $n\in\Z$, the triangle inequality, estimates 
(\ref{split1}), (\ref{split}) and (\ref{split3}) complete the proof of Lemma \ref{claim}.
\end{proof}
Next we apply Lemma \ref{claim} with $\Delta _{N_{j}L_{j}}(u_{j})$, $j=1,2$ in the place of $u_{j}$ 
and it follows from the definition of $X^{s,b}$ that, for any $b>1/2$ and any $\delta >0$,
\begin{equation}\label{vtora-pista}
\big\|\prod_{j=1}^{2}\Delta _{N_{j}L_{j}}(u_{j})\big\|_{L^{2}(\R\times M)}
\leq 
C_{b,\delta}N_{1}^{\frac{1}{2}+\delta}(L_1L_2)^{b}
\prod_{j=1}^{2}
\|\Delta _{N_{j}L_{j}}(u_{j})\|_{L^{2}(\R\times M)}\, .
\end{equation}
It is now clear that the proof of Proposition \ref{bolen} can be completed from~\eqref{parva-pista} and~\eqref{vtora-pista} 
by a suitable H\"older inequality
\end{proof}
Let us now turn to the proof of (\ref{parvo-red}). 
Set  $\partial=\frac{\partial}{\partial z}$ and $\overline{\partial}=\frac{\partial}{\partial \overline{z}}$.
Thanks to (\ref{gauge}) and using that $F(0)=0$, we obtain that the function
$$
F(u)-(\partial F)(0)\, u-(\overline{\partial}F)(0)\,\overline{u}
$$
is vanishing at least of order $3$ at the origin. 
Therefore, in order to prove (\ref{parvo-red}) , it suffices to prove 
\begin{equation}\label{parvo-red-bis}
\|F(u)\|_{X^{s,-b'}(\R\times M)}\leq C\Big(\|u\|_{X^{1,b}(\R\times M)}^{2}+\|u\|_{X^{1,b}(\R\times M)}^{p}\Big)
\|u\|_{X^{s,b}(\R\times M)}\
\end{equation}
assuming that $F(u)$ is vanishing to order $3$ in zero. We can write
$$
F(u)=\sum_{N_1}
\big[
F(S_{N_1}(u))-F(S_{N_1/2}(u))
\big],
$$
where the sum is taken over all dyadic values of $N_1$ (recall that $S_{1/2}(u)=0$).
We have for $z,w\in\C$,
$$
F(z)-F(w)=(z-w)\int_{0}^{1}\partial F(tz+(1-t)w)dt+(\overline{z}-\overline{w})\int_{0}^{1}\overline{\partial}F(tz+(1-t)w)dt\, .
$$
Therefore
\begin{multline*}
F(S_{N_1}(u))-F(S_{N_1/2}(u))\\
=
\Delta _{N_1}(u)G_{1}(\Delta _{N_1}(u),S_{N_1/2}(u))+\overline{\Delta _{N_1}(u)}G_{2}(\Delta _{N_1}(u),S_{N_1/2}(u)),
\end{multline*}
with
$$
G_{1}(z_1,z_2)=\int_{0}^{1}\partial F(tz_1+z_2)dt,\quad
G_{2}(z_1,z_2)=\int_{0}^{1}\overline{\partial} F(tz_1+z_2)dt.
$$
We have thus the splitting
$
F(u)=F_{1}(u)+F_{2}(u),
$
where
$$
F_{1}(u)=\sum_{N_1}\Delta _{N_1}(u)\,G_{1}(\Delta _{N_1}(u),S_{N_1/2}(u))\, .
$$ 
Thanks to the growth assumption on $V(z)$, we have the bound
$$
|G_{j}(z_1,z_2)|\leq C(1+|z_1|+|z_2|)^{\alpha-1},\quad j=1,2 .
$$
We will provide a bound only for $F_{1}(u)$. The analysis for $F_{2}(u)$ is exactly the same.

We have for dyadic integers $N_1$, $N_2$
\begin{equation*}
\Delta _{N_1}\,\Delta _{N_2}=
\begin{cases}
\Delta _{N_1} & \text{if $N_1=N_2$,}
\\
0 &\text{otherwise.}
\end{cases}
\end{equation*}
and hence we can write that for dyadic integers $N_1$, $N_2$,
\begin{equation*}
G_{1}(\Delta _{N_1}S_{N_2}(u),S_{N_{1}/2}S_{N_2}(u))-G_{1}(\Delta _{N_1}S_{N_{2}/2}(u),S_{N_{1}/2}S_{N_{2}/2}(u))
\end{equation*}
is equal to
\begin{equation*}
\begin{cases}
G_{1}(0,S_{N_2}(u))-G_{1}(0,S_{N_{2}/2}(u)) & \text{ if $2N_{2}\leq N_{1}$,}  \\
G_{1}(\Delta _{N_1}(u),S_{N_{1}/2}(u))-G_{1}(0, S_{N_{1}/2}(u)) & \text{ if $N_{2}=N_{1}$,}  \\
G_{1}(\Delta _{N_1}(u),S_{N_{1}/2}(u))-G_{1}(\Delta _{N_1}(u),S_{N_{1}/2}(u))=0 & \text{ if $N_{2}\geq 2N_1$}.
\end{cases}
\end{equation*}
Using the vanishing property of $F$ at the origin allows us to write
\begin{multline*}
G_{1}(\Delta _{N_1}(u),S_{N_1/2}(u))=
\sum_{N_2\,:\,N_{2}\leq N_{1}}\Delta _{N_2}(u)H_{1}^{N_2}(\Delta _{N_2}(u),S_{N_{2}/2}(u))
+
\\
+
\sum_{N_{2}\,:\,N_2\leq N_{1}}\overline{\Delta _{N_2}(u)}H_{2}^{N_2}(\Delta _{N_2}(u),S_{N_{2}/2}(u)),
\end{multline*}
with
\begin{equation*}
H_{1}^{N_2}(a,b)=
\begin{cases}
\int_{0}^{1}\partial_{2}G_{1}(0,ta+b)dt & \text{ if $2N_{2}\leq N_1$}, \\
\int_{0}^{1}\partial_{1}G_{1}(ta,b)dt & \text{ if $N_2=N_1$},
\end{cases}
\end{equation*}
where $(\partial_{1},\partial_{2})$ are the derivatives of $G_1$ with respect to the first and the second arguments respectively.
Moreover
\begin{equation*}
H_{2}^{N_2}(a,b)=
\begin{cases}
\int_{0}^{1}\overline{\partial_{2}}G_{1}(0,ta+b)dt & \text{ if $2N_{2}\leq N_1$}, \\
\int_{0}^{1}\overline{\partial_{1}}G_{1}(ta,b)dt & \text{ if $N_2=N_1$}.
\end{cases}
\end{equation*}
Notice that
$$
|H_{j}^{N_2}(a,b)|\leq C(1+|a|+|b|)^{\max(\alpha-2,0)},\quad j=1,2 .
$$
We can write
\begin{multline*}
F_{1}(u)=\sum_{N_2\leq N_1}\Delta _{N_1}(u)\Delta _{N_2}(u)H_{1}^{N_2}(\Delta _{N_2}(u),S_{N_{2}/2}(u))+
\\
+\sum_{N_2\leq N_1}\Delta _{N_1}(u)\overline{\Delta _{N_2}(u)}H_{2}^{N_2}(\Delta _{N_2}(u),S_{N_{2}/2}(u))
:=F_{11}(u)+F_{12}(u)\, .
\end{multline*}
We will provide a bound only for $F_{11}(u)$. The analysis for $F_{12}(u)$ is exactly the same.

Similarly to the analysis for $G_{1}(\Delta _{N_1}(u),S_{N_1/2}(u))$, using once again
the vanishing property of $F$ at the origin, allows us to expand $H_{1}^{N_2}$ as follows
\begin{multline*}
H_{1}^{N_2}(\Delta _{N_2}(u),S_{N_{2}/2}(u))=
\sum_{N_3\,:\,N_{3}\leq N_{2}}\Delta _{N_3}(u)H_{11}^{N_2,N_3}(\Delta _{N_3}(u),S_{N_{3}/2}(u))
+
\\
+
\sum_{N_{3}\,:\,N_3\leq N_{2}}\overline{\Delta _{N_3}(u)}H_{12}^{N_2,N_3}(\Delta _{N_3}(u),S_{N_{3}/2}(u)),
\end{multline*}
where, due to the growth assumptions on $V$, $H_{1j}^{N_2,N_3}(a,b)$ satisfies
\begin{equation}\label{alpha-2}
|H_{1j}^{N_2,N_3}(a,b)|\leq C(1+|a|+|b|)^{\max(\alpha-3,0)},\quad j=1,2 .
\end{equation}
Of course we can write explicit formulas for $H_{1j}^{N_2,N_3}(a,b)$ as we did for $H_{1}^{N_2}(a,b)$
but it will not be needed in the sequel. The only information for $H_{1j}^{N_2,N_3}(a,b)$ that we will use is the
bound (\ref{alpha-2}).
Now, we can write
\begin{multline*}
F_{11}(u)=\sum_{N_3\leq N_2\leq N_1}
\Delta _{N_1}(u)\Delta _{N_2}(u)\Delta _{N_3}(u)
H_{11}^{N_2,N_3}(\Delta _{N_3}(u),S_{N_{3}/2}(u))+
\\
+\sum_{N_3\leq N_2\leq N_1}\Delta _{N_1}(u)\Delta _{N_2}(u)
\overline{\Delta _{N_3}(u)}H_{12}^{N_2,N_3}(\Delta _{N_3}(u),S_{N_{3}/2}(u))
:=F_{111}(u)+F_{112}(u)\, .
\end{multline*}
We will provide a bound only for $F_{111}(u)$. The analysis for $F_{112}(u)$ is exactly the same.
Notice that
\begin{equation}\label{more}
\Delta _{N}=\sum_{L}\Delta _{NL},
\end{equation}
where the sum is taken over all dyadic values of $L$. For $w\in X^{-s,b'}(\R\times M)$, we set
$$
I:=\sum_{\genfrac{}{}{0pt}{}{L_{0},L_{1},L_2,L_3, N_0}{N_3\leq N_2\leq N_1}}\int_{\R\times M}\Delta _{N_0 L_0}(w)
\prod_{j=1}^{3}\Delta _{N_j L_j}(u)H_{11}^{N_2,N_3}(\Delta _{N_3}(u),S_{N_{3}/2}(u)),
$$
where the sum is taken over dyadic values of $N_j,L_j$, $j=0,1,2,3$. 
By duality, to prove (\ref{parvo-red-bis}) it suffices to establish the bound
$$
\left|I\right|\leq C\|w\|_{X^{-s,b'}(\R\times M)}
\Big(\|u\|_{X^{1,b}(\R\times M)}^{2}+\|u\|_{X^{1,b}(\R\times M)}^{p}\Big)\|u\|_{X^{s,b}(\R\times M)}\,.
$$
Set 
$$
I_{L_0,L_1,L_2,L_3}^{N_0,N_1,N_2,N_3}:=
\Bigl|\int_{\R\times M}\Delta _{N_0 L_0}(w)
\prod_{j=1}^{3}\Delta _{N_j L_j}(u)H_{11}^{N_2,N_3}(\Delta _{N_3}(u),S_{N_{3}/2}
(u))\Bigr|\, .
$$
We split $I$ as $|I|\leq I_{1}+I_{2}$,
where we define $I_1$ and $I_2$ to be the sums of the terms $I_{L_0,L_1,L_2,L_3}^{N_0,N_1,N_2,N_3}$ 
associated to indexes such that $N_0\leq \Lambda N_1$ and $ N_0> \Lambda N_1$ respectively, and $\Lambda>1$ 
is a large constant to be determined later.

We first evaluate $I_1$. Using Proposition \ref{bolen}, and the H\"older inequality, we get,
that for every $\varepsilon>0$ there exists $\beta<\frac{1}{2}$ such that,
\begin{multline*}
I_{L_0,L_1,L_2,L_3}^{N_0,N_1,N_2,N_3}\leq C_{\varepsilon}
(N_2 N_3)^{\frac{1}{2}+\varepsilon}(L_0 L_1 L_2 L_3)^{\beta}
\|H_{11}^{N_2,N_3}(\Delta _{N_3}(u),S_{N_{3}/2}(u))\|_{L^{\infty}(\R\times M)}
\\
\|\Delta _{N_0 L_0}(w)\|_{L^{2}(\R\times M)}\prod_{j=1}^{3}\|\Delta _{N_j L_j}(u)\|_{L^2(\R\times M)}
\end{multline*}
Thanks to (\ref{alpha-2}) we can write
$$
\|H_{11}^{N_2,N_3}(\Delta _{N_3}(u),S_{N_{3}/2}(u))\|_{L^{\infty}}\leq
C\big(1+\|\Delta _{N_3}(u)\|_{L^{\infty}}+\|S_{N_{3}/2}(u)\|_{L^{\infty}}\big)^{\max(\alpha-3,0)}
$$
Using Lemma \ref{sob-co-g}, (\ref{harcterizatzia}) and the Cauchy-Schwarz inequality yield, for $b>1/2$,
\begin{multline*}
\|\Delta _{N_3}(u)\|_{L^{\infty}(\R \times M)}
\leq \sum_{L}\|\Delta _{N_3 L}(u)\|_{L^{\infty}(\R \times M)}
\leq
\\
\leq
\sum_{L}
CN_{3}^{\frac{3}{2}}L^{\frac{1}{2}}\|\Delta _{N_3 L}(u)\|_{L^{2}(\R \times M)}
\leq
CN_{3}^{\frac{1}{2}}\sum_{L}L^{\frac{1}{2}-b}L^{b}N_{3}
\|\Delta _{N_3 L}(u)\|_{L^{2}(\R \times M)}
\\
\leq CN_{3}^{\frac{1}{2}}\Big(\sum_{L}L^{1-2b}\Big)^{\frac{1}{2}}
\Big(\sum_{L}L^{2b}N_{3}^{2}\|\Delta _{N_3 L}(u)\|^{2}_{L^{2}(\R \times M)}\Big)^{\frac{1}{2}}
\leq 
CN_{3}^{\frac{1}{2}}\|u\|_{X^{1,b}}
\end{multline*}
We next estimate $\|S_{N_{3}/2}(u)\|_{L^{\infty}}$. 
\begin{multline*}
\|S_{N_{3}/2}(u)\|_{L^{\infty}(\R\times M)}
\leq 
\\
\leq
\sum_{N_{4}\,:\,N_{4}\leq N_{3}/2}
\|\Delta _{N_4}(u)\|_{L^{\infty}(\R\times M)}
\leq C
\sum_{N_{4}\,:\,N_{4}\leq N_{3}/2}\sum_{L}N_{4}^{\frac{3}{2}}L^{\frac{1}{2}}\|\Delta _{N_4 L}(u)\|_{L^{2}(\R \times M)}
\leq
\\
\leq
C\Big(
\sum_{N_{4}\,:\,N_{4}\leq N_{3}/2}
\sum_{L}
(N_{4}^{\frac{1}{2}})^{2}L^{1-2b}
\Big)^{\frac{1}{2}}
\Big(
\sum_{N_{4}\,:\,N_{4}\leq N_{3}/2}
\sum_{L}
L^{2b}\,N_{4}^{2}\,\|\Delta _{N_4 L}(u)\|_{L^{2}}^{2}
\Big)^{\frac{1}{2}}
\leq 
\\
\leq
CN_{3}^{\frac{1}{2}}\|u\|_{X^{1,b}},
\end{multline*}
provided $b>1/2$.
Using the last two estimates, we obtain the bound
\begin{equation}\label{l-infty}
\|H_{11}^{N_2,N_3}(\Delta _{N_3}(u),S_{N_{3}/2}(u))\|_{L^{\infty}(\R\times M)}\leq
1+C(N_{3}^{\frac{1}{2}}\|u\|_{X^{1,b}})^{\max(\alpha-3,0)}\, .
\end{equation}
With (\ref{l-infty}) in hand, we estimate $I_1$. 
Let us recall a discrete Schur lemma.
\begin{lemme}\label{ds}
For every $\Lambda>0$, every $s>0$ there exists $C>0$ such that
if $(c_{N_{0}})$ and $(d_{N_1})$ are two sequences of nonnegative numbers
indexed by the dyadic integers, then,
\begin{equation*}
\sum_{N_{0}\leq \Lambda N_{1}}\, \frac{N_{0}^{s}}{N_1^{s}}\,c_{N_{0}}\,d_{N_1}
\leq C\Big(\sum_{N_0}c_{N_{0}}^{2}\Big)^{\frac{1}{2}}\Big(\sum_{N_1}d_{N_1}^{2}\Big)^{\frac{1}{2}}\, .
\end{equation*}
\end{lemme}
\begin{proof}
Let us set
$$
K(N_0,N_1):=
\11_{N_{0}\leq \Lambda N_{1}}\,\frac{N_{0}^{s}}{N_1^{s}}\, .
$$
Summing geometric series imply that there exists $C>0$ such that
$$
\sup_{N_0}\sum_{N_1}K(N_0,N_1)+\sup_{N_1}\sum_{N_0}K(N_0,N_1)\leq C\, .
$$
Therefore the Schur lemma implies the boundedness on $l^{2}_{N_0}\times l^{2}_{N_1}$  of the bilinear form with 
kernel $K(N_0,N_1)$. This completes the proof of Lemma \ref{ds}.
\end{proof}
In estimation $I_1$, we first sum with respect to $L_0,L_1,N_0,N_1$. Writing 
$$
1= \frac{N_{0}^{s}}{N_1^{s}}\, N_{0}^{-s}\,N_{1}^{s},\quad
(L_0 L_1)^{\beta}=L_{0}^{b'}L_{1}^{b}L_{0}^{\beta-b'}L_{1}^{\beta-b},
$$
using Lemma \ref{ds} and  (\ref{l-infty}), after summing geometric series in $L_0$, $L_1$, we can write for $b>1/2$ and $1/2>b'>\beta$,
\begin{multline*}
I_{1}\leq
C_{\varepsilon }
\|u\|_{X^{s,b}(\R\times M)}\|w\|_{X^{-s,b'}(\R\times M)}
\Big(1+\|u\|_{X^{1,b}(\R\times M)}^{\max(\alpha-3,0)}\Big)
\\
\sum_{L_2,L_3}\sum_{N_3\leq N_2}
(N_2 N_3)^{\frac{1}{2}+\varepsilon}
(L_2 L_3)^{\beta}N_{3}^{\frac{\max(\alpha-3,0)}{2}}
\prod_{j=2}^{3}\|\Delta _{N_j L_j}(u)\|_{L^2(\R\times M)}\, .
\end{multline*}
Since $\alpha<5$ and $N_{3}\leq N_2$,we have, choosing $\varepsilon>0$ small enough,
$$
(N_2 N_3)^{\frac{1}{2}+\varepsilon}\, N_{3}^{\frac{\max(\alpha-3,0)}{2}}\leq N_2 N_3 (N_2 N_3)^{-\varepsilon}\, .
$$
Therefore, by summing geometric series in $N_2,N_3,L_2,L_3$, we get the bound
$$
I_{1}\leq
C\|w\|_{X^{-s,b'}(\R\times M)}\Big(\|u\|_{X^{1,b}(\R\times M)}^{2}+
\|u\|_{X^{1,b}(\R\times M)}^{\alpha-1}\Big)\|u\|_{X^{s,b}(\R\times M)}\,.
$$
It remains to estimate $I_2$. This is performed by using the following proposition and summing geometric series. 
\begin{proposition}\label{prop.rough}
Let $s\geq 1$. Then there exists $\Lambda>0$ , $b,b'$ satisfying~\eqref{restrictions}, $\gamma>0$ and $p, C$ such that for every
$w\in X^{-s,b'}(\R\times M)$, $u\in X^{1,b}(\R\times M)$, if $N_{0}$, $N_1$, $N_2$, $N_3$
satisfy 
$$ N_0\geq \Lambda N_1, \quad N_3\leq N_2\leq N_1$$
then
\begin{equation}\label{eq.passirough}
 I_{L_0,L_1,L_2,L_3}^{N_0,N_1,N_2,N_3}
\leq
C(N_{0}L_{0}L_1 L_2 L_3)^{-\gamma}
\|w\|_{X^{-s,b'}}
\Big(\|u\|_{X^{1,b}}^{2}+\|u\|_{X^{1,b}}^{p}\Big)\|u\|_{X^{s,b}}\,.
\end{equation}
\end{proposition}
\begin{proof}To prove Proposition~\ref{prop.rough}, we consider three regimes:
\subsubsection{Case 1: $N_0^{1-\delta}\leq N_3$, $\delta>0$ small enough}

We use Proposition \ref{bolen}, the H\"older inequality and (\ref{l-infty}). For every $\varepsilon>0$, 
there exists $\beta<\frac{1}{2}$ such that,
\begin{multline*}
I_{L_0,L_1,L_2,L_3}^{N_0,N_1,N_2,N_3}\leq C_{\varepsilon}
(N_2 N_3)^{\frac{1}{2}+\varepsilon}(L_0 L_1 L_2 L_3)^{\beta}
N_{3}^{\frac{\max(\alpha-3,0)}{2}}
\|\Delta _{N_0 L_0}(w)\|_{L^{2}(\R\times M)}
\\
\prod_{j=1}^{3}\|\Delta _{N_j L_j}(u)\|_{L^2(\R\times M)}
\Big(1+\|u\|_{X^{1,b}(\R\times M)}^{\max(\alpha-3,0)}\Big).
\end{multline*}
Therefore
\begin{multline*}
I_{L_0,L_1,L_2,L_3}^{N_0,N_1,N_2,N_3}\leq C_{\varepsilon}
\frac{(L_0 L_1 L_2 L_3)^{\beta}}{L_{0}^{b'}(L_1 L_2 L_3)^{b}}
\frac{N_{0}^{s}}{N_{1}^{s}(N_2 N_3)}
(N_2 N_3)^{\frac{1}{2}+\varepsilon}
\\
N_{3}^{\frac{\max(\alpha-3,0)}{2}}
\|w\|_{X^{-s,b'}(\R\times M)}\Big(\|u\|_{X^{1,b}(\R\times M)}^{2}+
\|u\|_{X^{1,b}(\R\times M)}^{\alpha-1}\Big)\|u\|_{X^{s,b}(\R\times M)}\, .
\end{multline*}
Since $\alpha<5$, we observe that there exist $\varepsilon>0$, $\delta>0$ and $\gamma>0$ such that 
$$
\frac{N_{0}^{s}}{N_{1}^{s}(N_2 N_3)}
(N_2 N_3)^{\frac{1}{2}+\varepsilon}
N_{3}^{\frac{\max(\alpha-3,0)}{2}}
\leq
N_0^{-\gamma}\, .
$$
The parameter $\varepsilon>0$ being fixed, we choose $\beta$ as imposed by Proposition \ref{bolen}.
Finally we chose $(b,b')\in\R^2$ satisfying (\ref{restrictions}) such that $b'>\beta$.
\subsubsection{Case 2 : $N_3\leq N_0^{\frac 1 2}$}
 We start with a rough bound for $H_{11}^{N_2,N_3}$. 
By a repetitive use of Leibniz rule and the Sobolev embeddings, we obtain the following statement. 
\begin{lemme}\label{rough}
There exists $A>0$ such that for every coordinate patch 
$$
\kappa\,:\, U\subset \R^3\longmapsto M,
$$
every
$\gamma\in\N^3$, there exists $C_{\gamma}>0$ such that for every $u\in H^1(M)$, 
$$
\|\partial_{x}^{\gamma}H_{11}^{N_2,N_3}(\Delta _{N_3}(u(\kappa(x))),S_{N_{3}/2}(u(\kappa(x)))\|_{L^{\infty}(U)}
\leq C_{\gamma}N_{3}^{\frac{3|\gamma|} 2 +A}
\Big(1+\|u\|_{H^1(M)}^{|\gamma| +A}\Big)\, .
$$
\end{lemme}
We next state a bound for products of eigenfunctions.
\begin{lemme}\label{rough-bis}
Let $\frac 2 3 >\delta>0$.
There exists $\Lambda>0$ such that if 
\begin{equation}\label{ass-N}
N_{0}\geq \Lambda N_1,\quad N_3\leq N_{2}\leq N_{1},\quad N_{3}\leq N_{0}^{\frac 2 3 -\delta},
\end{equation}
then for every $\gamma>0$ there exists $C$ and $p$ such that for every $u,w\in L^{2}(M)$,
\begin{equation*}
\left|
\int_{M}
H_{11}^{N_2,N_3}(\Delta _{N_3}(u),S_{N_{3}/2}(u))\,
P_{k_{0}}w
\prod_{j=1}^{3}P_{k_j}u
\right|
\leq \frac{C}{N_{0}^{\gamma}}\|u\|_{L^{2}}^{3}\|w\|_{L^{2}}\Big(1+\|u\|_{H^1(M)}^{p}\Big),
\end{equation*}
provided $\langle\lambda_{k_{j}}\rangle^{\frac{1}{2}}
\in [N_{j},2N_{j}]$, $j=0,1,2,3$.
\end{lemme}
\begin{proof}
A similar argument already appeared in Lemma 2.6 of our previous paper \cite{BGT'03}. 
The new point here is the presence of $H_{11}^{N_2,N_3}$.
Working in local coordinates, due to Lemma \ref{sogge}, we can  substitute $P_{k_0}w$ with the oscillatory integral
\begin{equation}\label{oscil-1}
\int e^{i\lambda_{k_{0}}^{\frac{1}{2}}\varphi(x,y_0)}
a_{0}(x,y_0,\lambda_{k_{0}}^{\frac{1}{2}})w(y_0)dy_0\, .
\end{equation}
Indeed the remainder term can be estimated thanks to the Sobolev embeddings and Lemma \ref{rough}.
We consider three cases.
\begin{itemize}
\item{Case 1.}\label{case111}
Suppose first that $N_{1}\leq N_{0}^{1-\delta}$. 
Using Lemma \ref{pak}, we integrate by parts in the variable $x$ 
by means of the oscillating factor 
$$
e^{i\lambda_{k_{0}}^{\frac{1}{2}}\varphi(x,y_0)},
$$
and after $q$ integrations,
we gain a factor $N_{0}^{-q}$. On the other hand, due to Lemma \ref{rough}, the assumption $N_{1}\leq N_{0}^{1-\delta}$
and the Sobolev inequality, we obtain that the derivation of the amplitude is causing at most a factor $N_{0}^{q(1-\delta)+A}$. 
By taking $q\gg 1$, this completes the proof in the case $N_{1}\leq N_{0}^{1-\delta}$.

\item{Case 2.}\label{case222}
Suppose next that $N_{1}\geq N_{0}^{1-\delta}$ but $N_{2}\leq N_{0}^{1-\delta}$. In this case we can substitute 
$P_{k_0}w$ with (\ref{oscil-1}) and $P_{k_1}u$ with 
\begin{equation*}
\int e^{i\lambda_{k_{1}}^{\frac{1}{2}}\varphi(x,y_1)}
a_{1}(x,y_1,\lambda_{k_{1}}^{\frac{1}{2}})u(y_1)dy_1\, .
\end{equation*}
Indeed in the considered case the remainders in the approximation for $P_{k_0}w$ and $P_{k_1}u$
given by Lemma \ref{sogge} are both ${\cal O}(N_{0}^{-\infty})$ as operators from $L^2$ to the Sobolev spaces.
Thanks to Lemma \ref{pak}, if we take $\Lambda\gg 1$, we can again integrate by parts in $x$ 
with the slightly modified oscillatory factor
$$
e^{i\lambda_{k_{0}}^{\frac{1}{2}} \Phi(x,y_0,y_1)},
$$
where
$$
\Phi(x,y_0,y_1)=\varphi(x,y_0)+\lambda_{k_{0}}^{-\frac{1}{2}}\lambda_{k_{1}}^{\frac{1}{2}}  
\varphi(x,y_1)\, .
$$
\item{Case 3.}\label{case.333}
Suppose finally that $N_{1}\geq N_{0}^{1-\delta}$ and $N_{2}\geq N_{0}^{1-\delta}$. Then we can substitute 
$P_{k_0}w$, $P_{k_1}u$ and $P_{k_2}u$ with the corresponding oscillatory integrals and we can then argue as in case~2.
\end{itemize}
This completes the proof of Lemma \ref{rough-bis}.
\end{proof}
Lemma \ref{rough-bis} (with $\delta = \frac 1 6$) is now used to prove Proposition~\ref{prop.rough} for space time functions in this regime.

Define $\Pi_{k,L}$ as follows
$$
\Pi_{k,L}(u):=
\frac{1}{2\pi}\int_{L\leq\langle\tau+\lambda_k\rangle\leq 2L}\widehat{P_k u}(\tau)e^{it\tau}d\tau\, .
$$
Further we set
$$
\Lambda(N_0,N_1,N_2,N_3):=
\{(k_0,k_1,k_2,k_3)\,:\, N_j\leq \langle\lambda_{k_{j}}\rangle^{\frac{1}{2}}\leq 2N_j, \, j=0,1,2,3\}\, .
$$
Since
$$
\Delta _{NL}=\sum_{k\,:\,N\leq \langle\lambda_{k}\rangle^{\frac{1}{2}}< 2N} \Pi_{k,L}
$$
we get the bound\begin{multline*}
I_{L_0,L_1,L_2,L_3}^{N_0,N_1,N_2,N_3}
\\
\leq
\sum_{\Lambda(N_0,N_1,N_2,N_3)}
\Bigl|
\int_{\R\times M}
H_{11}^{N_2,N_3}(\Delta _{N_3}(u),S_{N_{3}/2}(u))\
\Pi_{k_0,L_0}(w)
\prod_{j=1}^{3}\Pi_{k_j, L_j}(u)
\Bigr|.
\end{multline*}
Since $P_k\Pi_{k,L}=\Pi_{k,L}$, under the assumption (\ref{ass-N}) a use of Lemma \ref{rough-bis} yields,
\begin{multline*}
I_{L_0,L_1,L_2,L_3}^{N_0,N_1,N_2,N_3}
\leq
C_{\gamma}N_{0}^{-\gamma}
\sup_{t\in\R}\Big(1+\|u(t)\|_{H^1(M)}^{p}\Big)
\\
\sum_{\Lambda(N_0,N_1,N_2,N_3)}
\int_{-\infty}^{\infty}
\|\Pi_{k_0,L_0}w(t)\|_{L^{2}(M)}\prod_{j=1}^{3}\|\Pi_{k_{j},L_{j}}u(t)\|_{L^{2}(M)}dt\, .
\end{multline*}
For $b>1/2$, a use of (\ref{Sobolev}) and the H\"older inequality implies that 
\begin{multline*}
I_{L_0,L_1,L_2,L_3}^{N_0,N_1,N_2,N_3}
\leq
C_{\gamma}N_{0}^{-\gamma}
\Big(1+\|u\|_{X^{1,b}(\R\times M)}^{p}\Big)
\\
\sum_{\Lambda(N_0,N_1,N_2,N_3)}
\|\Pi_{k_0,L_0}w\|_{L^{2}(\R\times M)}\|\Pi_{k_1,L_1}u\|_{L^{2}(\R\times M)}
\prod_{j=2}^{3}\|\Pi_{k_j,L_j}u\|_{L^{\infty}(\R\,;\, L^{2}(M))}\, .
\end{multline*}
Since $\Delta _{N_{j}L_{j}}\Pi_{k_j,L_j}=\Pi_{k_j,L_j}$, 
for $k_j$ such that 
$N_j\leq \langle\lambda_{k_{j}}\rangle^{\frac{1}{2}}< 2N_j$,
using (\ref{eq-bis}), we get
\begin{equation}\label{i-s}
\|\Pi_{k_0,L_0}w\|_{L^{2}(\R\times M)}\leq CN_{0}^{s}L_{0}^{-b'}\|w\|_{X^{-s,b'}(\R\times M)}\, ,
\end{equation}
\begin{equation}\label{ii-s}
\|\Pi_{k_1,L_1}u\|_{L^{2}(\R\times M)}\leq CN_{1}^{-s}L_{1}^{-b}\|u\|_{X^{s,b}(\R\times M)}
\leq CL_{1}^{-b}\|u\|_{X^{s,b}(\R\times M)}\, ,
\end{equation}
and using Lemma \ref{sob-co-g}, for $j=2,3$,
\begin{multline}\label{iii-s}
\|\Pi_{k_j,L_j}u\|_{L^{\infty}(\R\,;\, L^{2}(M))}\leq L_{j}^{\frac{1}{2}}\|\Pi_{k_j,L_j}u\|_{L^{2}(\R\times M)}
\\
\leq CL_{j}^{\frac{1}{2}-b}\|u\|_{X^{0,b}(\R\times M)}
\leq CL_{j}^{\frac{1}{2}-b}\|u\|_{X^{1,b}(\R\times M)}\, .
\end{multline}
Using a crude form of Weyl asymptotics, we get a bound
\begin{equation}\label{Weyl}
|\Lambda(N_0,N_1,N_2,N_3)|\leq C\big(\prod_{j=0}^{3}N_{j}\big)^{c}\, .
\end{equation}
Estimate~\eqref{eq.passirough} in this regime follows in view of (\ref{i-s}), (\ref{ii-s}), (\ref{iii-s}) and (\ref{Weyl}).
\subsubsection{Case 3: $N_0^{\frac 1 2} \leq N_3 \leq N_0^{1-\delta}$ where $\delta$ is the small number fixed in Case 1}
We shall denote by $\mathcal{O}(1)$ any quantity bounded by 
$$(L_0L_1L_2L_3)^{-\gamma }\, \| w\| _{X^{-s,b'}}(\| u\|
_{X^{1,b}}^2+\| u\|
_{X^{1,b}}^p)\| u\|
 _{X^{s,b}}\ $$
for some $\gamma >0, p \in \mathbb{N}$.
Let 
$$v= \Delta _{N_0L_0}(w)\, \Delta _{N_1L_1}(u)\Delta _{N_2L_2}(u)\Delta
_{N_3L_3}(u)\ .$$
\begin{lemme} There exists $\varphi \in C^\infty _0(\R \setminus \{
0\} )$ such that
$$\| (1-\varphi (N_0^{-2}{\mathbf \Delta }))v\| _{L^1(\R \times M)}\leq
\mathcal{O}(1)N_0^{-k}$$ for any $k$.
\end{lemme}
Indeed, working in local coordinates, according to~\cite[Proposition 2.1]{BGT1}, there exists $\chi\in C^\infty_0( \mathbb{R}^3 \setminus \{0\})$ such that for any $k$
$$\|v-\prod _{j=0}^3 \chi (N_j^{-1}D)(u_j)\|_{L^1(\mathbb{R}\times M)}\leq \mathcal{O}(1)N_0^{-k}$$
Therefore, modulo negligible terms, the spectrum of $v$ lies in a ring of size $N_0$ which proves the lemma. Next we take advantage of this spectral localization to perform integrations by parts: we have (for some function $\psi\in C^\infty_0( \mathbb{R}\setminus \{0\})$)
$$\varphi (N_0^{-2}{\mathbf \Delta })=N_0^{-2}{\mathbf \Delta }\circ \psi (N_0^{-2}{\mathbf \Delta })$$
and modulo negligible terms
$$I_{L_0, L_1, L_2, L_3} ^{N_0, N_1, N_2, N_3}=N_0^{-2}\left|\int _{\R \times M}v\, \psi (N_0^{-2}{\mathbf \Delta })[{\mathbf \Delta}
(H_{11}^{N_2, N_3}( \Delta _{N_3}(u), S_{N_3/2} (u)))]\, dt\, dx\right|$$
Applying Proposition~\ref{bolen} we obtain
$$I_{L_0, L_1, L_2, L_3} ^{N_0, N_1, N_2, N_3}\leq \mathcal{O}(1)\,  N_0^{s-2}{N_1}^{-s}\, N_2^{-\frac{1}{2}+\varepsilon }\,
N_3^{-\frac{1}{2}+\varepsilon }\, \| {\mathbf \Delta }(H_{11}^{N_2, N_3}( \Delta _{N_3}(u), S_{N_3/2} (u)))\|
_{L^\infty }\ ,$$
But, by Sobolev embedding, we have 
$$\| {\mathbf \Delta }(H_{11}^{N_2, N_3}( \Delta _{N_3}(u), S_{N_3/2} (u)))\|
_{L^\infty }\leq C N_3^3 (1+ \|u\|_{X^{1,b}}^2)$$
and thus, since $N_3\leq N_2\leq N_1$ and $N_3 \leq N_0^{1- \delta}$,
\begin{multline}I_{L_0, L_1, L_2, L_3} ^{N_0, N_1, N_2, N_3}
\leq \mathcal{O}(1)\,  N_0^{s-2}{N_1}^{-s}\, N_2^{-\frac{1}{2}+\varepsilon }\,
N_3^{\frac{5}{2}+\varepsilon }\\
\leq \mathcal{O}(1)\,  N_0^{(s-2)}\, 
N_3^{2-s+2\varepsilon } \leq N_0^{(s-2) \delta +2 \varepsilon}. 
\end{multline}
and we can choose $\varepsilon>0$ small enough such that $(s-2)\delta+2\varepsilon <0$.

This completes the proof of Proposition~\ref{prop.rough} (and thus of (\ref{parvo-red})) in the case 
of a three dimensional compact manifold satisfying Proposition~\ref{str-bl-s3}.\end{proof}
\begin{remarque}\label{zabelejka}
Let us notice that the estimate  (\ref{parvo-red}) holds for any sub-quintic nonlinearity, not necessarily satisfying the
gauge condition (\ref{gauge}). We used (\ref{gauge}) in
the reduction to $F$ vanishing of order three at zero performed above because it simplifies a bit the analysis.
More precisely for an arbitrary $F$ in the expansions of $G_{1}(\Delta _{N_1}(u),S_{N_1/2}(u))$ and 
$H_{1}^{N_2}(\Delta _{N_2}(u),S_{N_{2}/2}(u))$ above one should add a constant.
This would force one to analyze quadratic nonlinearities separately which can be done with our methods. 
\end{remarque}
Thanks to the multilinear nature of our arguments,
the proof of (\ref{vtoro-red}) is essentially the same as for (\ref{parvo-red}).
Indeed for suitable $F_1$, $F_2$ one writes 
$$
F(u)-F(v)=(u-v)F_1(u,v)+(\overline{u}-\overline{v})F_2(u,v).
$$
Then we expand
$$
u-v=\sum_{N_1}\Delta _{N_1}(u-v)
$$
and for $j=1,2$,
$$
F_{j}(u,v)=\sum_{N_2}\big[F_{j}(S_{N_2}(u), S_{N_2}(v))-F_{j}(S_{N_{2}/2}(u),S_{N_2/2}(v))\big].
$$
One then further expand the difference and after a duality argument the proof of (\ref{vtoro-red}) is reduced
to a bound for a $4$-linear expression multiplied with a factor similar to $H_{11}^{N_2,N_3}(\Delta _{N_3}(u),S_{N_{3}/2}(u))$ 
appeared in the proof of (\ref{parvo-red}). We omit the details.
\section{Trilinear Strichartz estimates and applications}
\label{sec.5}
In this section we prove Theorem~\ref{thm4} for $M=S^2_{\rho}\times S^1$ with the standard metric. 
\subsection{Trilinear Strichartz estimates on $M=S^2_{\rho}\times S^1$}
We do not know whether Proposition \ref{str-bl-s3} holds in this case. 
Instead, we shall prove a trilinear Strichartz-type estimate. Let us first introduce some notation.
As usual we identify $S^1$ with $\R/(2\pi\Z)$. The eigenfunctions of ${\mathbf \Delta}$ in the considered case are
$$
\lambda_{m,n}=m^{2}+\kappa (n^{2}+n)  ,\quad m\geq0,\quad n\geq 0,\,\quad \kappa = \frac 1 {\rho^2}\,.
$$ 
Let us denote by $\Pi_{n}$ the spectral projector on spherical harmonics of degree $n\geq 0$ on $S^2_{\rho}$.
For $f(\omega,\theta)\in L^{2}(S^{2}_{\rho}\times S^{1})$, we set
$$
\Theta_{m}f(\omega):=\frac 1 { 2 \pi}\int_{0}^{2\pi}f(\omega,\theta)\, e^{-im\theta}d\theta\, .
$$
The crucial estimate is the following.
\begin{proposition}\label{str-tl-s21}
For every interval $I\subset\R$, every $\varepsilon>0$ there exists a constant $C$ such that for every 
$N_1\geq N_2\geq N_3\geq 1$, every $f_1,f_2,f_3\in L^{2}(M)$,
\begin{equation*}
\big\|\prod_{j=1}^{3}e^{it{\mathbf \Delta}}(\Delta _{N_j}f_{j})\big\|_{L^{2}(I\times M)}
\leq C N_{3}^{\frac{5}{4}}N_{2}^{\frac{3}{4}+\varepsilon}\prod_{j=1}^{3}\|\Delta _{N_j}f_j\|_{L^{2}(M)}\, .
\end{equation*}
\end{proposition}
\begin{proof}
By a time translation we can suppose that $I=[0,T]$.
Since $\kappa$ is not necessarily integer, we can not employ the argument of
Proposition \ref{str-bl-s3} which reduces the analysis to the case $I=[0,2\pi]$. 
We shall instead use the following lemma, already used in a similar context in \cite{Bo4}.
\begin{lemme}\label{cska}
Let $\Lambda$ be a countable set of real numbers. Then for every $T>0$ there exists $C_T$ such that for every sequence $(a_{\lambda})$
indexed by $\Lambda$ one has
$$
\big\|\sum_{\lambda\in\Lambda}a_{\lambda}\, e^{i\lambda t}\big\|_{L^{2}(0,T)}
\leq C_{T}\Big(\sum_{l\in\Z}\Big(\sum_{\lambda\,:\,|\lambda-l|\leq 1/2}|a_{\lambda}|\Big)^{2}\Big)^{\frac{1}{2}}\, .
$$
\end{lemme}
\begin{proof}
Let $\psi_{T}\in C_{0}^{\infty}(\R)$ be such that $\psi_{T}=1$ on the interval $[0,T]$. Set
$$
f(t):=\sum_{\lambda\in\Lambda}\psi_{T}(t)\,a_{\lambda}\, e^{i\lambda t}\, .
$$
Then
$$
\widehat{f}(\tau)=\sum_{\lambda\in\Lambda}\widehat{\psi_{T}}(\tau-\lambda)a_{\lambda}
$$
and the problem is to show that
$$
\|\widehat{f}\|_{L^2(\R)}\leq 
 C_{T}\Big(\sum_{l\in\Z}\Big(\sum_{\lambda\,:\,|\lambda-l|\leq 1/2}|a_{\lambda}|\Big)^{2}\Big)^{\frac{1}{2}}\, .
$$
Next, we write
$$
|\widehat{f}(\tau)|\leq \sum_{l\in\Z}\,\,\,\sum_{\lambda\,:\,|\lambda-l|\leq 1/2}
|\widehat{\psi_{T}}(\tau-\lambda)|\, |a_{\lambda}|
\leq
\sum_{l\in\Z} \, K(l,\tau)h(l),
$$
where
$$
h(l)=\sum_{\lambda\,:\,|\lambda-l|\leq 1/2}|a_{\lambda}|\, ,
\qquad
K(l,\tau)=\sup_{\lambda\,:\, |\lambda-l|\leq 1/2}|\widehat{\psi_{T}}(\tau-\lambda)|\, .
$$
It is clear that $|\lambda-l|\leq 1/2$ implies
$$
\frac{1}{1+|\tau-\lambda|}\leq \frac{C}{1+|\tau-l|}
$$
and therefore, using that $\psi_{T}\in C_{0}^{\infty}(\R)$, we deduce that for every $N\in\N$ there exists $C_{T,N}$ such that
$$
|K(l,\tau)|\leq \frac{C_{T,N}}{(1+|\tau-l|)^{N}}\,\, .
$$
A use of the Schur lemma completes the proof of Lemma~\ref{cska}.
\end{proof}
We expand
\begin{multline*}
\Big(\prod_{j=1}^{3}e^{it{\mathbf \Delta}}(\Delta _{N_j}f_{j})\Big)(\omega,\theta)
=
\\
=
\sum
e^{-i(\lambda_{m_1,n_1}+\lambda_{m_2,n_2}+\lambda_{m_3,n_3})t}
e^{i(m_1+m_2+m_3)\theta}
\prod_{j=1}^{3}(\Pi_{n_{j}}\Theta_{m_{j}}f_{j})(\omega)\, ,
\end{multline*}
where the sum is taken over $(m_{j},n_{j})$, $j=1,2,3$ such that
$N_{j}\leq \langle\lambda_{m_{j},n_{j}}\rangle^{\frac{1}{2}}< 2N_j$.
Using the Parseval identity with respect to $\theta$ and Lemma~\ref{cska}, we obtain
\begin{multline*}
\big\|\prod_{j=1}^{3}e^{it{\mathbf \Delta}}(\Delta _{N_j}f_{j})\big\|^{2}_{L^{2}([0,T]\times M)}
\\
\leq
C_{T}\,\,
\sum_{(l,\xi)\in\Z^{2}}
\Big\|\displaystyle\sum_{\stackrel{|l-\lambda_{m_1,n_1}-\lambda_{m_2,n_2}-\lambda_{m_3,n_3}|\leq 1/2}{\xi=m_1+m_2+m_3}}
\,\,\,
\prod_{j=1}^{3}\Bigl|\Pi_{n_{j}}\Theta_{m_{j}}f_{j}\Bigr| \Big\|_{L^{2}(S^2_{\rho})}^{2}\,,
\end{multline*}
where the summation over $(m_1,m_2,m_3,n_1,n_2,n_3)$ is restricted to $(m_j,n_j)$ such that
$N_{j}\leq  \langle\lambda_{m_{j},n_{j}}\rangle^{\frac{1}{2}}< 2N_j$, $j=1,2,3$.
Applying the triangle inequality for the $L^{2}(S^2_{\rho})$ norm, 
the Cauchy-Schwarz inequality in the summation over $(m_1,m_2,m_3,n_1,n_2,n_3),$  
and the trilinear estimate (\ref{thm2-2}) of Theorem \ref{thm2} yields
\begin{multline*}
\big\|\prod_{j=1}^{3}e^{it{\mathbf \Delta}}(\Delta _{N_j}f_{j})\big\|^{2}_{L^{2}([0,T]\times M)}
\\
\leq
\sum_{(l,\xi)\in\Z^{2}}
\displaystyle\sum_{\stackrel{|l-\lambda_{m_1,n_1}-\lambda_{m_2,n_2}-\lambda_{m_3,n_3}|\leq 1/2}{\xi=m_1+m_2+m_3}}
|\Lambda(l,\xi)|
(N_2 N_3)^{\frac{1}{2}}\prod_{j=1}^{3}\|\Pi_{n_{j}}\Theta_{m_{j}}f_{j}\|_{L^{2}(S^2_{\rho})}^{2}
\\
\leq
(N_2 N_3)^{\frac{1}{2}}
\sup_{(l,\xi)\in\Z^{2}}|\Lambda(l,\xi)|\prod_{j=1}^{3}\|f_{j}\|_{L^{2}(M)}^{2}\, ,
\end{multline*}
where
\begin{multline*}
\Lambda(l,\xi):=\{(m_1, m_2, m_3, n_1, n_2,n_3)\in\N^{6}\, :\, 
\big|l-\sum_{j=1}^{3}(m_{j}^{2}+\kappa(n_j^2+n_j))\big|\leq 1/2,\quad
\\
\xi=m_1+m_2+m_3,\quad
N_{j}\leq  \langle\lambda_{m_{j},n_{j}}\rangle^{\frac{1}{2}}< 2N_j,
\quad
j=1,2,3
\}\, .
\end{multline*}
It remains to bound the size of $\Lambda(l,\xi)$.
The number of possible $(m_3,n_3)$ is bounded by $CN_3^2$. The number of possible $m_2$ is bounded by $CN_2$.
Thus the number of possible $(m_2,m_3,n_3)$ is bounded by $CN_{2}N_{3}^{2}$.
Let us now {\bf fix} a possible triple  $(m_2,m_3,n_3)$. Our goal is evaluate the number of possible
$(m_1,n_1,n_2)$ such that $(m_1, m_2, m_3, n_1, n_2,n_3)\in \Lambda(l,\xi)$. In view of the imposed restrictions, we can
eliminate $m_1$ by concluding that $(n_1,n_2)$ should satisfy
$$
\left|l-(\xi-m_2-m_3)^{2}-m_{2}^{2}-m_{3}^{2}-\kappa\big[n_{1}^{2}+n_{2}^{2}+n_{3}^{2}+n_1+n_2+n_3\big]\right|\leq \frac{1}{2}
$$
or equivalently
\begin{equation}\label{broene}
\left|(2n_1+1)^{2}+(2n_2+1)^{2}-R\right|\leq \frac{2}{\kappa},
\end{equation}
where 
$$
R=-4(n_3^2+n_3)+2+\frac{4}{\kappa}\big[l-(\xi-m_2-m_3)^{2}-m_{2}^{2}-m_{3}^{2}\big]\,.
$$
Using Lemma~\ref{le3.2}, uniformly with respect to $R$, 
the number of integer solutions $(n_1,n_2)\in [0, CN_1]\times [0, CN_2]$ of the inequality (\ref{broene}) is bounded
by $C_{\varepsilon}N_{2}^{\varepsilon}$ which implies the estimate
$$
|\Lambda(l,\xi)|\leq C_{\varepsilon}N_{3}^{2}N_{2}^{1+\varepsilon}\, .
$$
The proof of Proposition \ref{str-tl-s21} is now completed.
\end{proof}
\subsection{Using trilinear Strichartz estimates}
From now on we simply assume that $M$ is a three dimensional compact manifold satisfying Proposition~\ref{str-tl-s21}.
Proceeding as in~\cite[Section 3.2]{BGT'03} one can show, for instance, that the product of any  
Zoll surface with $S^1$ 
has this property. As a consequence it can be remarked that in fact Theorem~\ref{thm1} holds for any such manifold. 
 
For our purpose in this section, we will first use the following weaker form of Proposition~\ref{str-bl-s3} 
which is a consequence of Proposition \ref{str-tl-s21}. 
\begin{proposition}\label{str-bl-s21}
For every interval $I\subset\R$, every $\varepsilon>0$ there exists a constant $C$ such that 
for every $N_1,N_2\geq 1$, every $f_1,f_2\in L^{2}(M)$,
\begin{equation*}
\big\|\prod_{j=1}^{2}e^{it{\mathbf \Delta}}(\Delta _{N_j}f_{j})\big\|_{L^{2}(I\times M)}
\leq C (\min(N_1,N_2))^{\frac{3}{4}+\varepsilon}\prod_{j=1}^{2}\|\Delta _{N_j}f_j\|_{L^{2}(M)}\, .
\end{equation*}
\end{proposition}
\begin{proof}
It suffices to apply Proposition \ref{str-tl-s21} with $f_{3}=1$.
\end{proof}
Proposition \ref{str-tl-s21} and Proposition \ref{str-bl-s21} now imply the following statement.
\begin{proposition}\label{bolen-bis}
For every $\varepsilon>0$ there exist $\beta<\frac{1}{2}$ and $C>0$ such that for every $N_1\geq N_2\geq N_3\geq 1$,
$L_1,L_2, L_3\geq 1$, 
every $u_1,u_2,u_3\in L^{2}(\R\times M)$,
\begin{equation}\label{lille-1}
\big\|\prod_{j=1}^{2}\Delta _{N_{j}L_{j}}(u_{j})\big\|_{L^{2}}
\leq C(L_1 L_2 )^{\beta}
N_{2}^{\frac{3}{4}+\varepsilon}\prod_{j=1}^{2}\|\Delta _{N_j L_j}(u_j)\|_{L^{2}}
\end{equation}
and
\begin{equation}\label{lille-2}
\big\|\prod_{j=1}^{3}\Delta _{N_{j}L_{j}}(u_{j})\big\|_{L^{2}}
\leq C(L_1 L_2 L_3 )^{\beta}
N_{3}^{\frac{5}{4}+\varepsilon}
N_{2}^{\frac{3}{4}+\varepsilon}
\prod_{j=1}^{3}\|\Delta _{N_j L_j}(u_j)\|_{L^{2}}\, .
\end{equation}
\end{proposition}
\begin{proof}
One can show that  Proposition \ref{str-bl-s21} implies (\ref{lille-1}) exactly as we did in the proof of Proposition \ref{bolen}.
The proof of (\ref{lille-2}) follows similar lines. First, using  Lemma~\ref{sob-co-g} and the H\"older inequality we get
\begin{multline}\label{parva-pista-bis}
\big\|\prod_{j=1}^{3}\Delta _{N_{j}L_{j}}(u_{j})\big\|_{L^{2}}
\leq 
\|\Delta _{N_{1}L_{1}}(u_{1})\|_{L^{6}(\R\,;\,L^{2}(M))}
\prod_{j=2}^{3}\|\Delta _{N_{j}L_{j}}(u_{j})\|_{L^{6}(\R\,;\,L^{\infty}(M))}
\\
\leq
C\, (N_{2}N_{3})^{\frac{3}{2}}(L_1L_2L_3)^{\frac{1}{3}}
\prod_{j=1}^{3}\|\Delta _{N_j L_j}(u_j)\|_{L^{2}}\, .
\end{multline}
Next, exactly as in the proof of  Proposition \ref{bolen}, we obtain that for every unit interval $I\subset \R$, every
$b>1/2$, every $\delta>0$ there exists $C_{b,\delta}$ such that for every $N_{1}\geq N_{2}\geq N_{3}\geq 1$, every
$u_{1},u_{2},u_3\in X^{0,b}(\R\times M)$,
\begin{equation}\label{claim-bis}
\big\|\prod_{j=1}^{3}\Delta _{N_{j}}(u_{j})\big\|_{L^{2}(I\times M)}
\leq C_{b,\delta}\, N_2^{\frac{3}{4}+\delta}N_3^{\frac{5}{4}}\prod_{j=1}^{3}\|\Delta _{N_j}(u_j)\|_{X^{0,b}} \, .
\end{equation}
Using the  partition of unity (\ref{split}), we get the bound
\begin{multline}\label{vtora-pista-bis}
\big\|\prod_{j=1}^{3}\Delta _{N_{j}L_{j}}(u_{j})\big\|_{L^{2}(\mathbb{R}\times M)}\\
\leq 
C_{b,\delta}N_2^{\frac{3}{4}+\delta}N_3^{\frac{5}{4}}(L_1L_2L_{3})^{b}
\prod_{j=1}^{3}\|\Delta _{N_{j}L_{j}}(u_{j})\|_{L^{2}(\mathbb{R}\times M)}\, .
\end{multline}
Finally, a suitable interpolation between 
(\ref{parva-pista-bis}) and  (\ref{vtora-pista-bis}) completes the proof of Proposition \ref{bolen-bis}.
\end{proof}
Let us now turn to the proof of Theorem \ref{thm4} in the case $M=S^{2}_\rho\times S^{1}$ 
(or more generally any manifold satisfying Proposition~\ref{str-tl-s21}). We can again suppose that $F(u)$ is vanishing
at least at order three at zero. We expand $F(u)$ as we did in section~\ref{sec.4} and are led to estimating terms of the form
$$
I_{L_0,L_1,L_2,L_3}^{N_0,N_1,N_2,N_3}:=
\Bigl|\int_{\R\times M}\Delta _{N_0 L_0}(w)
\prod_{j=1}^{3}\Delta _{N_j L_j}(u)H_{11}^{N_2,N_3}(\Delta _{N_3}(u),S_{N_{3}/2}(u))\Bigr|\, .
$$
As in section~\ref{sec.4}, we consider two cases and denote by $I_1$ the contribution corresponding to $N_0\leq \Lambda N_1$, where $\Lambda\gg 1$ is a large constant and by $I_2$ the contribution corresponding to $N_0> \Lambda N_1$. To study $I_1$, we even make one more expansion of the terms
$H_{1j}^{N_2,N_3}$ and it results that estimate (\ref{parvo-red}) is a consequence
of the bounds
\begin{equation}\label{J}
J\leq C\|w\|_{X^{-s,b'}(\R\times M)}\|u\|_{X^{1,b}(\R\times M)}\|u\|_{X^{s,b}(\R\times M)}^{2}
\end{equation}
and
\begin{equation}\label{I1}
I_1\leq C\|w\|_{X^{-s,b'}(\R\times M)}
\Big(\|u\|_{X^{1,b}(\R\times M)}^{3}+\|u\|_{X^{1,b}(\R\times M)}^{\alpha-1}\Big)\|u\|_{X^{s,b}(\R\times M)}\,,
\end{equation}
\begin{equation}\label{I2}
I_2\leq C\|w\|_{X^{-s,b'}(\R\times M)}
\Big(\|u\|_{X^{1,b}(\R\times M)}^{2}+\|u\|_{X^{1,b}(\R\times M)}^{p}\Big)\|u\|_{X^{s,b}(\R\times M)}\,,
\end{equation}
where
$$
J=\sum_{L_{0},L_1,L_2,L_3}\sum_{N_0}\sum_{N_{3}\leq N_{2}\leq N_{1}}\Bigl|\int_{\R\times M}\Delta _{N_0 L_0}(w)
\prod_{j=1}^{3}\Delta _{N_j L_j}(u)\Bigr|
$$
and
\begin{multline*}
I_1=
\sum_{L_{0},L_1,L_2,L_3,L_4}\,
\sum_{N_0\leq \Lambda N_1}\,\,\sum_{N_{4}\leq N_{3}\leq N_{2}\leq N_{1}}
\\
\Bigl|\int_{\R\times M}\Delta _{N_0 L_0}(w)\Big(\prod_{j=1}^{4}\Delta _{N_j L_j}(u)\Big)
H^{N_{2},N_{3}, N_4}(\Delta _{N_4}(u),S_{N_{4}/2}(u))\Bigr|\, ,
\end{multline*}
\begin{multline*}
I_2=
\sum_{L_{0},L_1,L_2,L_3}\,
\sum_{N_0>\Lambda N_1}\,\,\sum_{ N_{3}\leq N_{2}\leq N_{1}}
\\
\Bigl|\int_{\R\times M}\Delta _{N_0 L_0}(w)\Big(\prod_{j=1}^{3}\Delta _{N_j L_j}(u)\Big)
H^{N_{2},N_{3}}(\Delta _{N_3}(u),S_{N_{3}/2}(u))\Bigr|\, ,
\end{multline*}
with sums taken over the dyadic values of $N_{j}$ and $L_{j}$, $j=0,1,2,3,4$. Moreover
$H^{N_{2},N_{3},N_{4}}(a,b)$ enjoys the bound
$$
|H^{N_{2},N_{3},N_{4}}(a,b)|\leq C(1+|a|+|b|)^{\max(\alpha-4,0)}\, .
$$
In addition for $4<\alpha<5$, we can further expand $H^{N_{2},N_{3},N_{4}}$ and we can get the bound
\begin{equation}\label{posledno}
\left|H^{N_{2},N_{3},N_{4}}(\Delta _{N_4}(u),S_{N_{4}/2}(u))\right|^{\frac{1}{\alpha-4}}
\leq
C\sum_{N_{5}\,:\, N_{5}\leq N_{4}}|\Delta _{N_5}(u)|\, .
\end{equation}
The proof of (\ref{J}) is a consequence of the bilinear estimate (\ref{lille-1}). More precisely, using~(\ref{lille-1}) and
the H\"older inequality, we obtain that for every $\varepsilon>0$ there exists $\beta<1/2$ such that
$$
\left|
\int\Delta _{N_0 L_0}(w)\prod_{j=1}^{3}\Delta _{N_j L_j}(u)
\right|
\leq
C(N_{2}N_{3})^{\frac{3}{4}+\varepsilon}
L_{0}^{\beta}\|\Delta _{N_0 L_0}(w)\|_{L^2}
\prod_{j=1}^{3}L_{j}^{\beta}\|\Delta _{N_j L_j}(u)\|_{L^2}\, .
$$
Since for $\varepsilon<1/4$ we have $\frac{3}{4}+\varepsilon<1$, we can complete the proof of (\ref{J}) as we did in section~\ref{sec.4}. 
A similar argument  (using both (\ref{lille-1}) and (\ref{lille-2})) is valid for~\eqref{I1}, if $\alpha \leq 4$.

To prove~\eqref{I1} if $4<\alpha <5$,
we use Proposition \ref{bolen-bis} in its full strength. Set
$$
I_{L_0L_1L_2L_3L_4}^{N_0N_1N_2N_3N_4}:=\Bigl| \int_{\R\times M}\Delta _{N_0 L_0}(w)\Big(\prod_{j=1}^{4}\Delta _{N_j L_j}(u)\Big)
H^{N_{2},N_{3},N_{4}}(\Delta _{N_4}(u),S_{N_{4}/2}(u))\Bigr|\, .
$$
In order to estimate $I_{L_0L_1L_2L_3L_4}^{N_0N_1N_2N_3N_4}$, we use the following form of H\"older's inequality.
\begin{equation}\label{eq.holder}
\forall \gamma \in ]0,1], \qquad \left|\int_{\R\times M}\, fg\right|
\leq\Big(\int_{\R\times M}|f|\Big)^{1-\gamma}
\Big(\int_{\R\times M}|f|
|g|^{\frac{1}{\gamma}}\Big)^{\gamma}\,.
\end{equation}
Since $4< \alpha <5$,  $\gamma= \alpha -4\in ]0,1[$. Using~\eqref{eq.holder}, we can write,
\begin{equation}\label{inter}
I_{L_0L_1L_2L_3L_4}^{N_0N_1N_2N_3N_4}\leq [J_{L_0L_1L_2L_3L_4}^{N_0N_1N_2N_3N_4}]^{1-\gamma}
[K_{L_0L_1L_2L_3L_4}^{N_0N_1N_2N_3N_4}]^{\gamma}\, ,
\end{equation}
where
$$
J_{L_0L_1L_2L_3L_4}^{N_0N_1N_2N_3N_4}=\int_{\R\times M}
\Big|\Delta _{N_0 L_0}(w)\Big(\prod_{j=1}^{4}\Delta _{N_j L_j}(u)\Big)\Big|.
$$
Thanks to (\ref{posledno}), the second factor $K_{L_0L_1L_2L_3L_4}^{N_0N_1N_2N_3N_4}$ in (\ref{inter})  enjoys the bound
$$
K_{L_0L_1L_2L_3L_4}^{N_0N_1N_2N_3N_4}\leq 
C\sum_{N_{5}\,:\, N_{5}\leq N_{4}}
\int_{\R\times M}
\Big|\Delta _{N_0 L_0}(w)\Big(\prod_{j=1}^{4}\Delta _{N_j L_j}(u)\Big)\Big|\,
|\Delta _{N_5}(u)|\, .
$$
Let us now bound $J_{L_0L_1L_2L_3L_4}^{N_0N_1N_2N_3N_4}$. Using  H\"older inequality and
Proposition \ref{bolen-bis} (both (\ref{lille-1}) and (\ref{lille-2})),
we obtain that for every $\varepsilon>0$ there exists $\beta<1/2$ such that
$$
J_{L_0L_1L_2L_3L_4}^{N_0N_1N_2N_3N_4}\leq 
C(N_2 N_3)^{\frac{3}{4}+\varepsilon}N_{4}^{\frac{5}{4}+\varepsilon} L_{0}^{\beta}\|\Delta _{N_0 L_0}(w)\|_{L^2}
\prod_{j=1}^{4}L_{j}^{\beta}\|\Delta _{N_j L_j}(u)\|_{L^2}\, .
$$
Next we estimate $K_{L_0L_1L_2L_3L_4}^{N_0N_1N_2N_3N_4}$.
By writing $\Delta _{N_{5}}=\sum_{L_5}\Delta _{N_5 L_5}$, 
using H\"older's inequality and Proposition \ref{bolen-bis} (twice (\ref{lille-2})),
we obtain that for every $\varepsilon>0$ there exists $\beta<1/2$ such that
$K_{L_0L_1L_2L_3L_4}^{N_0N_1N_2N_3N_4}$ is bounded by
$$
C\sum_{N_{5}\,:\, N_{5}\leq N_{4}}\sum_{L_{5}}
(N_2 N_3)^{\frac{3}{4}+\varepsilon}(N_{4}N_{5})^{\frac{5}{4}+\varepsilon} 
L_{0}^{\beta}\|\Delta _{N_0 L_0}(w)\|_{L^2}
\prod_{j=1}^{5}L_{j}^{\beta}\|\Delta _{N_j L_j}(u)\|_{L^2}\, .
$$
Writing $N_{5}^{\frac{5}{4}+\varepsilon}=N_{5}N_{5}^{\frac{1}{4}+\varepsilon}$,
using the Cauchy-Schwarz inequality, we get for $b>\beta$, 
\begin{multline*}
\sum_{N_{5}\,:\, N_{5}\leq N_{4}}\sum_{L_{5}}N_{5}^{\frac{5}{4}+\varepsilon}
L_{5}^{\beta}\|\Delta _{N_5 L_5}(u)\|_{L^2}
\leq
\\
\leq
\Big(
\sum_{N_{5}\,:\, N_{5}\leq N_{4}}\sum_{L_{5}}\,\, \Big[L_{5}^{\beta-b}N_{5}^{\frac{1}{4}+\varepsilon}\Big]^{2}
\Big)^{\frac{1}{2}}\|u\|_{X^{1,b}}
\leq
CN_{4}^{\frac{1}{4}+\varepsilon}\|u\|_{X^{1,b}}\, .
\end{multline*}
Therefore, we have the estimate,
\begin{multline*}
K_{L_0L_1L_2L_3L_4}^{N_0N_1N_2N_3N_4}\leq 
C(N_2 N_3)^{\frac{3}{4}+\varepsilon}N_{4}^{\frac{5}{4}+\varepsilon} 
L_{0}^{\beta}\|\Delta _{N_0 L_0}(w)\|_{L^2}
\\
\Big(\prod_{j=1}^{4}L_{j}^{\beta}\|\Delta _{N_j L_j}(u)\|_{L^2}\Big)
N_{4}^{\frac{1}{4}+\varepsilon}\|u\|_{X^{1,b}}
\, .
\end{multline*}
Coming back to (\ref{inter}), we obtain the following estimate
\begin{multline}\label{crucial}
I_{L_0L_1L_2L_3L_4}^{N_0N_1N_2N_3N_4}\leq 
C\frac{N_0^{s}}{N_1^{s}}
\frac{(N_2 N_3)^{\frac{3}{4}+\varepsilon}N_{4}^{\frac{5}{4}+\varepsilon} N_{4}^{\gamma(\frac{1}{4}+\varepsilon)}}
{N_{2} N_{3} N_{4}}L_{0}^{\beta-b'}(L_1L_2L_3L_4)^{\beta-b}
\\
\big(N_{0}^{-s}L_{0}^{b'}\|\Delta _{N_0 L_0}(w)\|_{L^2}\big)\big(N_{1}^{s}L_{1}^{b}\|\Delta _{N_1 L_1}(u)\|_{L^2}\big)
\Big(\prod_{j=2}^{4}N_{j}L_{j}^{b}\|\Delta _{N_j L_j}(u)\|_{L^2}\Big)
\|u\|_{X^{1,b}}^{\gamma}\, .
\end{multline}
Let us take $\varepsilon>0$ such that $(\frac{3}{2}+2\varepsilon)+(\frac{5}{4}+\varepsilon)+\gamma(\frac{1}{4}+\varepsilon)<3$
or equivalently,
$$
0<\varepsilon<\frac{1-\gamma}{4(3+\gamma)}=\frac{5-\alpha}{4(\alpha-1)} \, .
$$
Note that a proper choice of $\varepsilon$ is possible thanks to the sub critical assumption $\alpha<5$. 
Therefore there exists $\theta>0$ such that for $N_{4}\leq N_{3}\leq N_{2}$,
\begin{equation}\label{glad}
\frac{(N_2 N_3)^{\frac{3}{4}+\varepsilon}N_{4}^{\frac{5}{4}+\varepsilon} N_{4}^{\gamma(\frac{1}{4}+\varepsilon)}}
{N_{2} N_{3} N_{4}}\leq\frac{C}{(N_{2} N_{3} N_{4})^{\theta}}\, .
\end{equation}

Thanks to (\ref{glad}) and (\ref{crucial}), we obtain
\begin{multline}I_{L_0L_1L_2L_3L_4}^{N_0N_1N_2N_3N_4}\leq 
C\frac{N_0^{s}}{N_1^{s}}\frac 1 {(N_2 N_3 N_4)^\theta} L_{0}^{\beta-b'}(L_1L_2L_3L_4)^{\beta-b}
\\
\big(N_{0}^{-s}L_{0}^{b'}\|\Delta _{N_0 L_0}(w)\|_{L^2}\big)\big(N_{1}^{s}L_{1}^{b}\|\Delta _{N_1 L_1}(u)\|_{L^2}\big)
\Big(\prod_{j=2}^{4}N_{j}L_{j}^{b}\|\Delta _{N_j L_j}(u)\|_{L^2}\Big)
\|u\|_{X^{1,b}}^{\gamma}\, .
\end{multline}
for $I_{1}$ by summing geometric series in $L_0$, $L_1$, $L_2$, $L_3$, $L_4$, $N_2$, $N_3$, $N_4$ while the sum over $(N_0,N_1)$ is 
performed by invoking Lemma \ref{ds}. 
\par

We now turn to the proof of~\eqref{I2}. As in section~4, we shall denote by $\mathcal{O}(1)$ any quantity bounded (for some $\gamma >0, p \in \mathbb{N}$) by 
$$(L_0L_1L_2L_3)^{-\gamma }\, \| w\| _{X^{-s,b'}}(\| u\|
_{X^{1,b}}^2+\| u\|
_{X^{1,b}}^p)\| u\|
 _{X^{s,b}}. $$

We have three regimes: 
\begin{enumerate} 
\item{ $N_0^{1- \delta}\leq N_3 $, $\delta> 0$ small enough}, 
\item{$N_3\leq N_0^{\frac 1 2}$},
\item {$N_0^{\frac 1 2} \leq N_3 \leq N_0^{1- \delta}$}.
\end{enumerate}
In the first regime, we apply the same strategy as when $N_0\leq \Lambda N_1$. Indeed, in this regime, $N_0/N_1\leq N_0^{\delta}$ and we obtain (after expanding once more the non linear term) with $\eta>0$
\begin{equation*}I_{L_0L_1L_2L_3L_4}^{N_0N_1N_2N_3N_4}\leq 
\mathcal{O}(1) \frac{N_0^{s}}{N_1^{s}}\frac 1 {(N_2 N_3 N_4)^\theta} \leq \mathcal{O}(1) N_0^\delta N_3^{-\theta}\leq\mathcal{O}(1) N_0^{\delta- \theta(1- \delta)}\leq \mathcal{O}(1)N_0^{-\eta}
\end{equation*}
which gives the summability in $N_4\leq N_3\leq N_2\leq N_1\leq \Lambda^{-1} N_0$.

The second regime can be dealt with in the same way as in the previous section (in this regime, we gain arbitrary powers of $N_0^{-1}$). Finally we concentrate on the last regime. 
Let 
$$v= \Delta _{N_0L_0}(w)\, \Delta _{N_2L_2}(u)\Delta
_{N_3L_3}(u)\ .$$
As in section~4, we have
\begin{lemme} There exists $\varphi \in C^\infty _0(\R \setminus \{
0\} )$ such that
$$\| (1-\varphi (N_0^{-2}{\mathbf \Delta }))v\| _{L^1(\R \times M)}\leq
\mathcal{O}(1)N_0^{-k}$$ for any $k$.
\end{lemme}
Using that for some function $\Psi\in C^\infty_0$ we have 
$$\varphi(N_0^{-2} {\mathbf \Delta})= N_0^{-2} {\mathbf \Delta }\circ \Psi(N_0^{-2} {\mathbf \Delta})$$
we can integrate by parts in the integral defining $I_{L_0,L_1,L_2,L_3}^{N_0, N_1, N_2, N_3}$:
\begin{multline}\label{eq.final}
I_{L_0,L_1,L_2,L_3}^{N_0,N_1,N_2,N_3}
= N_0^{-2}\int_M \Psi(N_0^{-2} {\mathbf \Delta})(v) {\mathbf \Delta }\bigl(\Delta _{N_1L_1}(u) H_{11}^{N_2,N_3}( \Delta _{N_3}(u), S_{N_3}(u))\bigr)\\
=  N_0^{-2}\int_M \Psi(N_0^{-2} {\mathbf \Delta})(v) \times \Bigl[{\mathbf \Delta }\left(\Delta _{N_1L_1}(u)\right) H_{11}^{N_2,N_3}( \Delta _{N_3}(u), S_{N_3}(u)) \\
+ \nabla\left(\Delta _{N_1L_1}(u)\right) \cdot \nabla \left(H_{11}^{N_2,N_3}( \Delta _{N_3}(u), S_{N_3}(u))\right)\\
+ \Delta _{N_1L_1}(u){\mathbf \Delta }\left(H_{11}^{N_2,N_3}( \Delta _{N_3}(u), S_{N_3}(u))\right)\Bigr] 
\end{multline}
Now to estimate the first term in the right hand side of~\eqref{eq.final} we simply apply the strategy already used in the case $N_0 \leq \Lambda N_1$, the only difference being the additional factor $(N_1/N_0)^{2}$ (the factor $N_1^2$ coming from the action of the Laplace operator on $\Delta_{N_1 L_1} (u)$), which allows, since $s<2$, to exchange the roles of $N_0$ and $N_1$ and gain summability. Remark here that the additional operator $\Psi(N_0^{-2} {\mathbf \Delta})$ applied to $\Delta _{N_0L_0}(w)\, \Delta _{N_2L_2}(u)\Delta _{N_3L_3}(u)$ plays no role since it disappears when one takes $L^2$ norms. 

Next we estimate the two other terms. We expand the derivatives and observe that
\begin{gather*}
 | \nabla \left(H_{11}^{N_2,N_3}( \Delta _{N_3}(u), S_{N_3}(u))\right)|(x)  \leq C w_{N_3,1}(x)\\
\text{ with }
w_{N_3,1} (x) = \left( | \nabla\Delta _{N_3}(u)| +  | \nabla S_{N_3}(u)|\right)\left( |\Delta _{N_3}(u)| +  |S_{N_3}(u)|\right) (x),
\end{gather*}
\begin{gather*}
 | {\mathbf \Delta}  \left(H_{11}^{N_2,N_3}( \Delta _{N_3}(u), S_{N_3}(u))\right)|(x)  \leq w_{N_3, 2} (x) \\
\text{ with }w_{N_3, 2}(x) = \left( |{\mathbf \Delta}\Delta _{N_3}(u)| +  | {\mathbf \Delta} S_{N_3}(u)|\right)\left( | \Delta _{N_3}(u)| +  | S_{N_3}(u)|\right) (x)\\
+ \left( | \nabla\Delta _{N_3}(u)| +  | \nabla S_{N_3}(u)|\right)^2(x).
\end{gather*}
Next we use that for $\phi \in C^\infty_0(\mathbb{R})$, $\phi(N^{-2} \mathbf{\Delta})$ is an $N^{-1}$-semi-classical operator and consequently the gradient of a spectrally localized function is essentially spectrally localized. This allows to apply the trilinear estimate~\eqref{lille-2} to $v$ on one hand and to 
\begin{equation}\label{eq.mod}
\begin{aligned}|\nabla\left(\Delta _{N_1L_1}(u)\right) | \times w_{N_3, 1} \\
\text{ or }
 \Delta _{N_1L_1}(u)\times w_{N_3,2} 
\end{aligned}
\end{equation}
 on the other hand (remark that the moduli in~\eqref{eq.mod} do not spoil the estimate since we take $L^2$ norms).

 We obtain that the contribution of these two terms is bounded (for any $\varepsilon>0$) by
\begin{multline} C_\epsilon \mathcal{O}(1) N_0^{-2} \left( \frac{N_0} {N_1}\right)^s N_2^{-1/4+ \epsilon} N_3^{1/4+ \epsilon} \left( N_1 N_3+ N_3^2\right) N_3^{-1/4+\varepsilon}N_3^{1/4+ \epsilon}\\
\leq C_\epsilon \mathcal{O}(1)N_0^{s-2} N_1^{1-s} N_3^{4\varepsilon +1}\leq C_\epsilon \mathcal{O}(1) N_0^{s-2} N_3^{2-s + 4 \varepsilon}
\end{multline} 
Using that $\delta >0$ and $N_3 \leq N_0^{1- \delta}$, for $\varepsilon >0$ small enough, this term can be bounded by $N_0^{- \eta}, \eta >0$ giving the required summability in $\Lambda^{-1}N_0\geq  N_1 \geq N_2 \geq N_3$.
\begin{remarque}
A careful examination of the proof above shows that  Theorem~\ref{thm4} still holds for a three dimensional manifold $M$ satisfying
the more general trilinear Strichartz estimate,
\begin{multline*}
\exists\, a>0\, :\, \forall\,  T>0,\,\, \forall\, \varepsilon>0,\, \exists\, C>0\, :\, 
\forall\, N_{3}\leq N_{2}\leq N_{1},\, \forall\,\, f_1,f_2,f_3\in L^2(M)\, ,
\\
\big\|\prod_{j=1}^{3}e^{it{\mathbf \Delta}}(\Delta _{N_j}f_{j})\big\|_{L^{2}([0,T]\times M)}
\leq C N_{3}^{1+a}N_{2}^{1-a+\varepsilon}\prod_{j=1}^{3}\|\Delta_{N_j}f_j\|_{L^{2}(M)}\, .
\end{multline*}
\end{remarque} 
\appendix
\section{}
This appendix is devoted to the optimality of the assumption $\alpha<5$ in Theorem~\ref{thm1}.
Let us again consider a 3d-manifold $M$ endowed with a Riemannian metric $g$ and ${\mathbf \Delta}$ the Laplace-Beltrami 
operator acting on functions of $M$. We consider the following non-linear Schr\"odinger equation on $M$
\begin{equation}\label{eq.surcrit}
(i\partial_t+{\mathbf \Delta}) u = F(u), u|_{t=0}= u_0 \in H^1(M)
\end{equation}
where $F(z)= (1+ |z|^2)^ {\frac{\alpha-1}{2}}z$ and $\alpha>5$.

Let us fix $s>3/2$. Equation (\ref{eq.surcrit}) is well-posed for data in $H^s(M)$ by the energy method.
In particular, for every bounded set $B\subset H^{s}(M)$ there exists $T_{s}$ such that for every $u_0\in H^{s}(M)$
the Cauchy problem (\ref{eq.surcrit}) has a unique solution on the interval $[-T_s,T_s]$ in the class
$C([-T_s,T_s]\,;\,H^{s}(M))$. Moreover the flow map
$$
\Phi\,:\, u_0\longrightarrow u
$$
is continuous (and even Lipschitz continuous) from $B$, endowed with the $H^{s}(M)$ metric, to $C([-T_s,T_s]\,;\,H^{s}(M))$.
The next statement shows that $\Phi$ can not be extended as a continuous map on bounded sets of $H^1(M)$.
\begin{theoreme}\label{th.surc}
Let $B$ be a bounded set of $H^1(M)$. There is no $T>0$ such that the map $\Phi$ can be extended as a continuous map from
$B$ to $C([-T,T]\,;\,H^{1}(M))$.
\end{theoreme}
The result of Theorem \ref{th.surc} readily follows from the following statement.
\begin{theoreme}\label{thm6}
There exist a sequence $(t_n)_{n\in\N}$ of positive numbers tending to zero and a sequence $(u_{n}(t))_{n\gg 1}$ of $C^{\infty}(M)$
functions defined for $t\in [0,t_n]$, such that 
$$
(i\partial_t+{\mathbf \Delta}) u_n =(1+ |u_n|^2)^ {\frac{\alpha-1}{2}}u_n
$$
with
$$
\lim_{n\rightarrow \infty}\|u_{n}(0)\|_{H^1(M)}=0\, ,\text{ and }
\lim_{n\rightarrow \infty}\|u_{n}(t_n)\|_{H^1(M)}= \infty \,.
$$
\end{theoreme}
\begin{remarque}
The result of Theorem \ref{thm6}, in the particular case $M= \R^3$, endowed with the standard metric,
can be found in \cite{CCT}. We also refer to  \cite{CCT} for more ill-posedness results for NLS on $\R^d$, $d\geq 1$,
with power-like nonlinearities and data in $H^s$, $s>0$.  
\end{remarque}
\begin{remarque}
The proof of Theorem \ref{thm6} is strongly inspired by \cite{CCT}. The only observation we make here
is that the dilation arguments involved in the proof in  \cite{CCT} are not essential.
It is clear from the proof we present that the discussed phenomenon is completely local, i.e. the whole analysis
is close to a point of $M$ for very small times.
\end{remarque}
\begin{proof}
We work in a local coordinate patch around $0$ and consider as initial data the sequence $u_{n}(0)= \kappa_n n^{1/2} \varphi(n x),\quad n\gg 1$, 
where $\varphi$ is a fixed non negative smooth compactly supported $\kappa_n= \log^{- \delta}(n)$ with $\delta>0$ to be fixed later. Remark that 
$$
\|u_{n}(0)\|_{H^1(M)} \sim \kappa_n.
$$
Let us set $f(z):=(1+ |z|^2)^ {\frac{\alpha-1}{2}}$.
Then 
$$
v_{n}(t)=\kappa_n n^{1/2} \varphi(n x) e^{-it f(\kappa_n n^{1/2}\varphi(n x))}
$$ 
is the solution of the equation
\begin{equation}
i\partial_t  v_n = F(v_n),\quad v_n|_{t=0}= u_{n}(0)\, .
\end{equation}
Let us give a basic bound for $v_{n}(t)$.
\begin{lemme}\label{le.expl}
There exist $c>0$ and $C>0$ such that for any $t\geq 0$,
$$
\|\nabla_x v_n(t)\|_{L^2} \geq \kappa_n \Big(c\,t\, \kappa_n^ {\alpha-1} n^{(\alpha-1)/2}-C\Big)\, .
$$
\end{lemme}
\begin{proof}
The change of variable $y= nx$ gives
\begin{equation}
\begin{aligned}
\|\nabla_x v_n(t)\|_{L^2}&= 
\kappa_n \Big\|\nabla_y \big[ \varphi(y) e^{-it f(\kappa_n n^{1/2}\varphi(y))}\big]\Big\|_{L^2}
\\
&\geq \smash{\kappa_n \bigl(2t \kappa_n n^{1/2}\,\|\varphi(y) \partial_z f(\kappa_n n^{1/2}\varphi(y))\cdot \nabla_y \varphi(y) \|_{L^2} - 
\|\nabla_y \varphi\|_{L^2}\bigr)}
\\
&\geq c\,t\,\kappa ^{\alpha}_n n^{(\alpha-1)/2} - C\kappa_n 
\end{aligned}
\end{equation}
which implies the lemma.
\end{proof}
For functions $u$ on $M$, we define the quantity,
$$
E_n(u):= \Big[n^2 \|u\|^2_{L^2} + n^{-2} \|{\mathbf \Delta }u\|^2_{L^2}\Big]^{\frac{1}{2}}\, .
$$
The key point in the proof of Theorem~\ref{thm6} is the next statement.
\begin{lemme}\label{le.est}
The solution $u_n$ of \eqref{eq.surcrit} with initial data $u_{0}=\kappa_n n^{1/2} \varphi(n x)\in C^{\infty}(M)$
exists for $0\leq t\leq t_n$, with $t_n = \log^{1/8}(n)n^{-(\alpha-1)/2}$.
Moreover, there exist $\epsilon >0$ such that for $t\in[0,t_n]$,
$$
E_n(u_n(t)- v_n(t))\leq C n^{-\epsilon}\, .
$$
\end{lemme}
\begin{proof}
Since the initial data are in $H^{s}$, $s>3/2$, we know that $u_{n}(t)$ exist on a (non empty) maximal time interval $[0,\widetilde{t}_n[$.
Consequently, to prove Lemma~\ref{le.est}, we simply prove the {\em a priori} estimates which ensure, by a classical
bootstrap argument, both the existence and the control on $E_n(u_n(t)- v_n(t))$ for $t\in[0,t_n]\cap [0, \widetilde {t}_n[$.
Let us set $w_n:=u_n -v_n$. 
For the sake of conciseness, in the rest of the proof of Lemma~\ref{le.est}, we drop the subscript $n$ of $u_n$, $v_n$ and $w_n$.
The a priori estimates involved in the proof are simply energy inequalities in the equations
\begin{align*}
(i \partial_t + {\mathbf \Delta})w &= F(u) - F(v) - {\mathbf \Delta }v= \mathcal{O} ( 1+ |v|^ {\alpha-1} + |w|^{\alpha-1}) w - {\mathbf \Delta }v\\
(i\partial_t + {\mathbf \Delta}){\mathbf \Delta }w &={\mathbf \Delta }\left(F(u) - F(v)\right)- {\mathbf \Delta}^2 v
=- {\mathbf \Delta}^2 v+\Lambda,
\end{align*}
where
\begin{multline*}
\Lambda:=\mathcal{O} (1+ |w|^{\alpha-1} +|v|^{\alpha-1}) {\mathbf \Delta }w 
+ 
\\
+
\mathcal{O}( (1+ |w|^{\alpha-2} + |v|^{\alpha-2})(1+ |w| + |v|+|\nabla w| + |\nabla v|)) \nabla w 
+
\\
+ \mathcal{O}\left((1+ |w|^{\alpha -3} +|v|^{\alpha-3})( 1+ |\nabla v| ^2 + (|v|+ |w|) |\nabla ^2 v|) \right) w\, .
\end{multline*}
From the explicit formula for $v$, we have for $0\leq t \leq t_n$, $k=0,1,2,\dots$,
$$
\|\nabla ^k v\|_{L^\infty} \leq C n^{1/2+ k} \log^{k/8}(n)\, .
$$
According to the Gagliardo-Nirenberg inequalities,
\begin{equation}\label{GN}
\|f\|_{L^{\infty}}\leq C \|f\|^{3/4}_{H^{2}} \|f\|^{1/4}_{L^2}\leq C n^{1/2} E_n (f)
\end{equation}
we deduce
\begin{multline}\label{eq.est1}
n\,\|\mathcal{O} ( 1+ |v|^ {\alpha-1} + |w|^{\alpha-1}) w\|_{L^2}
\leq C\left(1+\|v\|^{\alpha-1}_{L^\infty}+\|w\|^{\alpha-1}_{L^\infty}\right)\,\, n\|w\|_{L^2} \leq
\\
\leq C\,n^{(\alpha-1)/2}(E_n (w)+E_n^{\alpha}(w))\, .
\end{multline}
To estimate $\Lambda$, we proceed similarly. More precisely, 
thanks to (\ref{GN}),
we estimate systematically the terms involving $v$ or $w$ below the $\mathcal{O}$ sign in $L^\infty$. 
The only term which cannot be estimated by invoking (\ref{GN}) is 
\begin{equation}\label{new}
\mathcal{O}( (1+ |w|^{\alpha-2} + |v|^{\alpha-2})(|\nabla w|)) \nabla w. 
\end{equation}
In order to evaluate (\ref{new}), we use the bound
$$
\|\nabla w\|_{L^4}\leq C \|\nabla w\|_{H^{3/4}} \leq C n^{3/4} E_n(w)
$$
and we obtain
\begin{equation}\label{eq.est3}
\|(\ref{new})\|_{L^2}
\leq C\,n^{(\alpha+1)/2}(E_n (w)+E_n^{\alpha}(w))\, .
\end{equation} 
We are therefore conducted to the following estimate for $\Lambda$
\begin{equation}\label{eq.est2}
n^{-1}\, \|\Lambda\|_{L^2}
\leq C\, n^{(\alpha-1)/2} \log^{1/4}(n)
(E_n (w)+E_n^{\alpha}(w))\, .
\end{equation} 
Next, thanks to the formula for $v$,  for $0\leq t \leq t_n$, we estimate the source terms,
\begin{equation}
\label{eq.est4} 
n\|{\mathbf \Delta }v\|_{L^2} + n^{-1}\|{\mathbf \Delta }^2 v\|_{L^2}\leq  Cn^2 \log^{1/2}(n)
\end{equation}
According to~\eqref{eq.est1},~\eqref{eq.est3} and~\eqref{eq.est4}, we obtain
\begin{equation*}
\frac d {dt} E^{2}_{n}(w) 
\leq  C\, n^{(\alpha-1)/2} \log^{1/4}(n)(E_n^{2} (w)+E_n^{\alpha+1}(w))+Cn^2 \log^{1/2}(n)E_n(w)  \, .
\end{equation*}
Suppose first that $E_{n}(w) \leq 1$ which is clearly the case at least for $t\ll 1$ since $w|_{t=0}=0$.
Notice that
$$
2n^2 \log^{1/2}(n)E_n(w)\leq n^{(\alpha-1)/2} \log^{1/4}(n)E_n^{2}(w)
+
\frac{n^4 \log(n)}{n^{(\alpha-1)/2} \log^{1/4}(n)}.
$$
Therefore
\begin{equation*}
\frac d {dt} 
\Big[
e^{-C\,t\, n^{(\alpha-1)/2} \log^{1/4}(n)}
E^{2}_{n}(w) 
\Big]
\leq  C\, 
n^{4-\frac{\alpha-1}{2}}\log^{3/4}(n)\,
e^{-C\,t\, n^{(\alpha-1)/2} \log^{1/4}(n)} \, .
\end{equation*}
Integrating the last inequality between $0$ and $t$ gives the estimate
\begin{equation*}
E_{n}(w) \leq C\,n^{2-\frac{\alpha-1}{2}} \log^{1/4}(n)\,e^{ C\,t\, n^{(\alpha-1)/2} \log^{1/4}(n)} \, .
\end{equation*}
For every $\gamma>0$ there exists $C_{\gamma}$ such that for $t\in [0,t_n]$, 
$$
C\,t\, n^{(\alpha-1)/2} \log^{1/4}(n)\leq C\,\log^{3/8}(n)\leq \gamma \log n +C_{\gamma}
$$ 
Since $\alpha >5$, by taking $\gamma>0$ small enough, we obtain that there exists $\varepsilon>0$ such that for $t\in [0,t_n]$,
we have
\begin{equation}
E_{n}(w)\leq C\,n^{-\varepsilon}\, .
\end{equation}
Finally the usual bootstrap argument allows to drop the assumption $E_{n}(w) \leq 1$. 
This completes the proof of Lemma~\ref{le.est}.
\end{proof}
It is clear that by interpolation, the quantity $E_n(u)$ controls uniformly with respect to $n$ the $H^1$ norm of $u$.
Consequently, it follows from Lemma~\ref{le.expl} and Lemma~\ref{le.est} that for $\delta<\frac{1}{8\alpha}$ and $n\gg 1$,
$$
\|u_{n}(t_{n})\|_{H^{1}}\geq 
C\log^{\frac{1}{8}-\alpha\delta}(n)\, .
$$
This completes the proof of Theorem~\ref{thm6}.
\end{proof}
\backmatter

\end{document}